\documentclass[10pt,a4paper,reqno]{amsart}
\usepackage{mathrsfs}
\usepackage{amssymb}
\usepackage[utf8]{inputenc}
\usepackage[english]{babel}
\usepackage{amsfonts}
\usepackage{amsmath}
\usepackage{braket, mathtools, commath, hyperref}
\usepackage{thmtools}
\usepackage{float}
\DeclareFontFamily{U}{mathx}{\hyphenchar\font45}
\DeclareFontShape{U}{mathx}{m}{n}{
      <5> <6> <7> <8> <9> <10>
      <10.95> <12> <14.4> <17.28> <20.74> <24.88>
      mathx10
      }{}
\DeclareSymbolFont{mathx}{U}{mathx}{m}{n}
\DeclareFontSubstitution{U}{mathx}{m}{n}
\DeclareMathAccent{\widecheck}{0}{mathx}{"71}
\DeclareMathAccent{\wideparen}{0}{mathx}{"75}

\usepackage{amsthm}
\usepackage{color}
\usepackage{enumitem}
\usepackage[nameinlink]{cleveref}

\makeatletter
\def\@setauthors{%
  \begingroup
  \centering
  {\@author \par} % Nombre del autor en negrita y grande
  \vskip 1em
  {\addresses \par} % Afiliación en cursiva
  \vskip 1em
  \endgroup
}
\def\@setaddresses{%
  \par\nobreak \begingroup
  \small
  \raggedright
  \def\author##1{}% Evitar que el nombre se repita
  \def\\{\par\nobreak}%
  \addvspace\bigskipamount
  \addresses
  \par\endgroup
}
\makeatother

\author{Manuel Ca\~nizares}
\address{\textbf{Manuel Ca\~nizares}\textsuperscript{1}\\{\small\itshape BCAM - Basque Center for Applied Mathematics \\ Avenida Mazarredo 14. 48009 Bilbao, Spain.\\}}
\email{{\small manuel.canizares@ricam.oeaw.ac.at}}
\date{\today}
\thanks{\textsuperscript{1}Current address: {\itshape RICAM - Johann Radon Institute for Computational and Applied Mathematics. Altenberg Str. 69. 4040 Linz, Austria.}}
\title[Boundary deformation for a Neumann problem]{Boundary deformation techniques for Neumann problems for the Helmholtz equation}
\newtheorem{theorem}{Theorem}[section]
\newtheorem{definition}[theorem]{Definition}
\newtheorem{lemma}[theorem]{Lemma}
\newtheorem{cor}[theorem]{Corollary}
\newtheorem{prop}[theorem]{Proposition}
\newtheorem*{thm:direct}{Theorem \ref{thm:1.1}}
\newtheorem*{thm:inverse}{Theorem \ref{thm:inverse}}

\theoremstyle{remark}

\newtheorem*{acknowledgements}{Acknowledgements}

\newcommand{\be}{\begin{equation}}
\newcommand{\ee}{\end{equation}}
\newcommand{\supp}{\mathrm{supp}\,}
\newcommand*{\defeq}{\mathrel{\vcenter{\baselineskip0.5ex \lineskiplimit0pt
			\hbox{\scriptsize.}\hbox{\scriptsize.}}}%
	=}

\newcommand{\dd}{\,\mathrm{d}}

\newcommand\restr[2]{{% we make the whole thing an ordinary symbol
		\left.\kern-\nulldelimiterspace % automatically resize the bar with \right
		#1 % the function
		\right|_{#2} % this is the delimiter
}}
\renewcommand{\norm}[1]{\lVert#1\rVert}

\newcommand\helm{(\Delta+\lambda-V)\,}

\newcommand{\Rd}{\mathbb{R}^n}

\newcommand{\half}{{1/2}}

\newcommand{\D}{\Omega}
\newcommand{\C}{\mathcal{C}}

\newcommand{\R}{\mathbb{R}}
\newcommand{\N}{\mathbb{N}}

\newcommand{\Oprime}{\Omega\setminus\mathcal{V}}
\newcommand{\Vcal}{\mathcal{V}}

\newcommand{\X}{\mathbb{X}}

\begin{document}

\begin{abstract}
  We adapt boundary deformation techniques to solve a Neumann problem for the Helmholtz equation with rough electric potentials in bounded domains. In particular, we study the dependance of Neumann eigenvalues of the perturbed Laplacian with respect to boundary deformation, and we illustrate how to find a domain in which the Neumann problem can be solved for any energy, if there is some freedom in the choice of the domain. This work is motivated by a Runge approximation result in the context of an inverse problem in point-source scattering with partial data.
\end{abstract}
\maketitle
\section{Introduction}

In this paper we will study a boundary value problem of Neumann type for the Helmholtz equation with a compactly-supported electric potential, of the form
\begin{equation}
  \label{eq:neumann}
  \begin{cases}
      (\Delta+\lambda-V)\,u=f_1 & \textrm{in }\Omega,\\
      \partial_\nu u=f_2 & \textrm{on }\partial\Omega,    
  \end{cases}
\end{equation}
where $\Omega$ is a bounded open domain in $\Rd$, $n\geq3$, containing the support of the potential $V$. Here $\lambda>0$ is the energy, and $f_1$ and $f_2$ are given functions. Note that this problem doesn't have a solution in general. The goal of this work, however, is to show that \eqref{eq:neumann} can always be solved if one has some freedom in the choice of the domain. Roughly, we aim to answer the following question:

\textit{Given $\lambda$ and $V$, can one find a domain $\Omega$ such that the problem \eqref{eq:neumann} has a unique solution $u$ for any $f_1\in L^2(\Omega)$ and $f_2\in L^2(\partial\Omega)$?}

We will show that the answer is affirmative for big enough values of the energy, under low-regularity assumptions on the potential. In particular, we will consider potentials of the form \begin{equation}\label{V}V=V^0+\gamma^s+\alpha \dd\sigma,\end{equation} where 
\begin{itemize}
\item $V^0\in L^{n/2} (\Rd;\mathbb{R}),$ and it's compactly supported,
\item $\dd\sigma$ denotes the surface measure of a compact hypersurface $\Gamma$ which is locally described by the graph of Lipschitz functions and $\alpha\in L^\infty(\Gamma;\mathbb{R})$, and
\item $\gamma^s$ is of the form 
\[ \gamma^s=\chi^2\, D^sg, \] for some $s<1$, $g\in L^\infty(\Rd;\mathbb{R})$, and $\chi\,\in \mathcal{C}_0^\infty(\Rd;[0,1])$ is a compactly supported cut-off function. Here $D^s$ denotes the fractional derivative, defined via Fourier transform as $\widehat{D^s f}(\xi)=|\xi|^s\hat{f}(\xi)$.
\end{itemize}

In a given domain $\Omega$, the Neumann problem above can be solved for every choice of $f_1$ and $f_2$ as long as $\lambda$ is not a Neumann eigenvalue (NEV) for the operator $-\Delta+V$ in $\Omega$.  We say that a number $\lambda\in\R$ is a NEV for $-\Delta+V$ in $\Omega$ if there exists $\phi\in H^1(\Omega)$ not identically zero solving the homogeneuos Neumann problem.
\begin{equation}\label{eq:homogeneousneumann}
    \begin{cases}
            \left(\Delta +\lambda- V\right)\phi= 0 & \textrm{in }\Omega,\\
            \partial_{\nu}\phi=0 & \textrm{on }\partial \Omega,
        \end{cases}
\end{equation}
and the space of such $\phi$ is its corresponding eigenspace. We will call $\lambda$ a simple NEV if the corresponding eigenspace has dimension 1. 

Whenever $\lambda$ is not a NEV, one can solve \eqref{eq:neumann} using the method of layer potentials. 
With this in mind, if we fix $\lambda>0$, our goal will be to find a domain $\Omega$ in which $\lambda$ is not a NEV. We will give a semi-constructive argument to find such a domain, with an approach based in boundary deformation techniques. Indeed, we will first choose an arbitrary domain and then show how to perturb it by means of a diffeomorphism. 

In case $\lambda$ is a simple NEV, we will use the techniques in \cite{Henry2005} and \cite{henrotvariation2005} to find a formula for its derivative with respect to the perturbation of the domain, and prove that we can choose such a perturbation so that this derivative doesn't vanish. In particular, for any $\C^2$ vector field $\X$ supported away from the potential, we may deform the domain $\Omega$ through a family of diffeomorphisms of the form $h_t=i_\Omega+t\X$, for a fixed vector field $\X$ and $t$ small enough. We will then show that  there is a unique differentiable function $t\mapsto \lambda(t)$ such that $\lambda(t)$ is an eigenvalue of $-\Delta+V$ in $h_t(\Omega)$, $\lambda(0)=\lambda$ and its derivative at $t=0$ is
\begin{equation*}
  \dot \lambda(0)=\int_{\partial \Omega} \left(|\nabla_{\partial\Omega}u|^2-\lambda u^2\right)\X\cdot \nu,
\end{equation*}
where $u$ is the normalized eigenvector of $\lambda$ and $\nu$ is the normal outward-pointing vector of $\Omega$. Then, it will be a matter of choosing $\X$ in a way that the integral above doesn't vanish

We will also show that the property of having simple NEVs is generic, this is, most of the domains will have simple NEVs. In particular, we will show that the set of perturbations that produce multiple eigenvalues is meager (of first category), in the appropiate space of perturbations. This notion will be defined later. The main ingredient here is a generalization of Smale-Sard theorem proved by Dan Henry in \cite{Henry2005}. This theorem controls the size of the critical set of a map between Banach manifolds. The idea is to characterize multiple eigenvalues as critical points of such a map and then use this theorem. 

\subsection{Application to a Runge approximation}

Our initial motivation to answer this question comes from an inverse problem of electric scattering with local near-field data. In this problem, which was studied by the author in \cite{Caizares2024}, one attempts to identify an electric potential of compact support by placing point sources and measuring time-harmonic waves in points close to the support of the potential. To be precise, one has access to the scattered wave $u_{sc}$, which is the solution to the equation
\begin{equation}\label{eq:scatt}
  \begin{cases}
          \left(\Delta +\lambda- V\right)u_{sc}(\centerdot\,,y) = V \Phi_\lambda(\centerdot-y) & \textrm{in }\mathbb{R}^n,\\
          u_{sc}(\centerdot\,,y)\textrm{ satisfying SRC}.
      \end{cases}
  \end{equation}
  Here $\Phi_\lambda$ is the fundamental solution to the free Helmoltz equation at energies $\lambda>0$, which solves the distributional problem
\begin{equation*}\label{eq:homo}
    \begin{cases}
            \left(\Delta +\lambda\right)\Phi_\lambda = \delta_0 & \textrm{in }\mathbb{R}^n,\\
            \Phi_\lambda\textrm{ satisfying SRC}.
        \end{cases}
\end{equation*}
  SRC stands for the Sommerfeld Radiation Condition, which is a condition of decay at infinity. More precisely, a function $u$ is said to satisfy SRC if
\begin{equation}\label{eq:src}
  \lim_{|x|\to\infty}|x|^{\frac{d-1}{2}}\left(\frac{x}{|x|}\cdot\nabla u(x)-i\lambda^{1/2}u(x)\right)=0
\end{equation}
uniformly in every direction. This condition is physically meaningful, since solutions that satisfy it represent waves that ``radiate energy at infinity'' \cite{schoteighty1992}.
Note that $u_{sc}(x,y)$ can be interpreted as placing a point source in $y$ and measuring the scattered wave at the point $x$

The existence of the scattering wave as solution of the problem \eqref{eq:scatt} was given by the author in \cite{Caizares2024} by adapting an argument by Caro and García \cite{caroscattering2020}, who studied a similar inverse problem with potentials of the form $V=V^0+\alpha\dd\sigma$. They defined a family of functional-analytic spaces in which one can exploit appropiate resolvent estimates from harmonic analysis \cite{agmonasymptotic1976,keniguniform1987,ruizlocalnodate}. These estimates give a decay in suitable norms in terms on $\lambda$, which in turn means that the scattering solution can be constructed for big enough energies. Throughout this work, we will need to assume that \eqref{eq:scatt} has an unique solution, so that an assumption of the form $\lambda>\lambda^V$ will be present, where $\lambda^V$ is such that the scattering solution exists.

In our inverse problem, the denomination ``local near-field data'' refers to the fact that the measurements are taken in arbitrary small sets of codimension $1$. In particular, we aim to prove that, if $V_1$ and $V_2$ are two potentials of the form \eqref{V}, then
\[\restr{u_{sc,1}(x,y)}{\Sigma_1\times\Sigma_2}=\restr{u_{sc,2}(x,y)}{\Sigma_1\times\Sigma_2}\implies V_1=V_2,\]
where $u_{sc,j}$ is the scattering solutions associated to $V_j$, and $\Sigma_j$ are two open sets of codimension 1 and of class $C^3$. To prove identifiability with this data, the author made use of a Runge approximation from single-layer potentials with densities supported in $\Sigma_j$. This relied on the existence of a unique solution to the Neumann problem \eqref{eq:neumann} in a domain $\Omega$ whose boundary contained $\Sigma_1$ and $\Sigma_2$. The set of NEVs is countable and has no accumulation points, which means that, by choosing any such domain arbitrarily, it is possible to prove uniqueness for \textit{most} energies, this is, for all $\lambda>\lambda^V$ except for those that turn out to be NEVs of $-\Delta+V$ in $\Omega$.

The assumption that $\lambda$ is not an eigenvalue has been taken in other problems of inverse scattering, such as in \cite{harrachmonotonicity2019}.
However, in the aforementioned inverse problem, the choice of the domain in which one performs the Runge approximation -and thus that in which \eqref{eq:neumann} has to be solved- is quite arbitrary. Therefore, it makes sense to fix any energy $\lambda$, and attempt to find a suitable domain $\Omega$.  

A similar problem was considered by Stefanov in \cite{Stefanov1990}, where he sought to avoid Dirichlet eigenvalues of the Helmholtz equation in the context of an inverse problem in electric scattering. His argument is relatively straightforward as a consequence of the fact that Dirichlet eigenvalues are strictly monotonically decreasing with respect to domain inclusion \cite{Leis1967}. However, the setting of NEVs proves to be more complex. For instance, this monotonicity already does not hold in general for NEVs of the Laplacian \cite{Funano2016}.
\subsection{The result}

The main theorem of this work is stated below. Note that the second part of the theorem is needed in the inverse problem to perform the Runge approximation from the measurement sets $\Sigma_j$.

\begin{restatable}{theorem}{nevthm}\label{thm:nev}
  Let $n\geq 3$ and $V$ be a potential of the form \eqref{V}, and fix $\lambda>\lambda^V$.  
  
    We can find a bounded open domain $\Omega$ of class $\C^3$ such that $\supp V\subset \Omega$, in which there exists a unique solution $u\in H^1(\Omega)$ to the problem  
    \[
      \begin{cases}
              \left(\Delta +\lambda- V\right)u= f_1 & \textrm{in }\Omega,\\
              \partial_{\nu}u=f_2 & \textrm{on }\partial \Omega,
          \end{cases}
    \]
    for any $f_1\in H^{-1}_0(\Omega)\cap L^2(\Omega\setminus\supp V)$ and $f_2\in L^2(\partial \Omega)$.

    Furthermore, let $\Sigma$ be a non-empty set of dimension $n-1$, open in the $(n-1)$-dimensional topology, separated from $\supp V$, and that can be expressed as the graph of $\C^3$ functions. There exists $\Sigma'\subset\Sigma$ relatively open and non-empty such that we can find such a domain $\Omega$ as above, and satisfying that $\Sigma'\subset\partial\Omega$. 
\end{restatable}

\subsection{Outline of the paper}
\begin{itemize}
  \item In Section \ref{sect:preliminaries} we give some preliminary results that will be of use throughout the paper. In Section \ref{sect:ucp}, we give a series of unique continuation properties for the operator $\Delta-\lambda-V$. In Section \ref{sect:direct} we recall the solution of the scattering problem \eqref{eq:scatt} as given in \cite{Caizares2024}, with the definition of Caro and Garcia's functional-analytical spaces. This will be relevant in the subsequent sections.
  \item In Section \ref{sect:layer}, we use the method of layer potentials to characterize the conditions in which the Neumann problem \eqref{eq:neumann} is solvable in terms of the eigenspace associated to $\lambda$. This will allow us to solve it in the case in which $\lambda$ is not a NEV and develop the perturbation of eigenvalues later on.
  \item In Section \ref{sect:calculus} we will lay out the approach of boundary perturbations, and study the differentiation of differential operators and boundary conditions with respect to these perturbations. We end the section finding a formula for the derivative of a simple eigenvalue in Section \ref{sect:perturbation}.
  \item In Section \ref{sect:nev}, we show that the simplicity of eigenvalues is generic in the set of suitable perturbations. We start by giving the statement of Henry's genericity theorem \cite{Henry2005}, and recalling the notions of meager set, Lindel\"of space and semi-Fredholm operator in Section \ref{sect:henry}. Then, in Section \ref{sect:genericity}, we use this theorem to prove this genericity property.
  \item Finally, we give the proof of Theorem \ref{thm:nev} in Section \ref{sect:thm}.
  \item In Appendix \ref{sect:banach}, we recall the notion of Banach manifolds and their tangent space, based on the exposition by Lang \cite{Lang1999}.
  % \item We will give some basic notions of differential geometry in Appendix \ref{sect:diffgeo}. 
\end{itemize}

\begin{acknowledgements}
This work was developed during a stay of the author in Universitè de Bordeaux, under the supervision of Sylvain Ervedoza. The author would like to thank Sylvain for his guidance and support. 

This work was carried out under the financial support of the MCIN/AEI under FPI fellowship PRE2019-091776 and the project PID2021-122156NB-I00, and also by the Basque Government through the BERC 2022-2025 program and by the Ministry of Science and Innovation: BCAM Severo Ochoa accreditation CEX2021-001142-S a/ MICIN / AEI / 10.13039/501100011033.
\end{acknowledgements}

\section{Preliminaries}\label{sect:preliminaries}
\subsection{Unique continuation properties}\label{sect:ucp}
We pull now a Carleman estimate from \cite{Caizares2024}, which was obtained by perturbating such an estimate that Caro and Rogers obtained in \cite{caroglobal2016} for the Laplacian. This estimate is done in a family of Bourgain-type spaces that were introduced by Caro and Rogers, inspired by the works of Haberman and Tataru \cite{habermanuniqueness2013,habermanuniqueness2015}.
For $s\in\mathbb{R}$ and $\zeta\in\mathbb{C}^n$ we define the inhomogeneous Bourgain-type space $X_\zeta^s$ as the space of distributions $u\in\mathscr{S}'(\Rd)$ such that $\widehat{u}\in L^2_\text{loc}(\Rd)$ and 
\[ \norm{u}_{X_\zeta^s}=\norm{(M|\Re(\zeta)|^2+M^{-1}|p_\zeta|^2)^{s/2}\,\widehat{u}}_{L^2}<\infty, \]
endowed with the norm $\norm{\centerdot}_{X_\zeta^s}$, with $M>1$, where $\Re$ denotes the real part and 
\begin{equation}\label{eq:pzeta} p_\zeta(\xi)=-|\xi|^2+2i\zeta\cdot\xi+\zeta\cdot\zeta. \end{equation} 
The Carleman estimate can be stated as follows:
\begin{prop}{\cite{Caizares2024}}
  \label{lem:carlemanest}
        Let $R_0>0$ such that $\supp V\subset B_{R_0}=\{x\in\Rd:|x|<R_0\}$. For $\zeta\in \mathbb{C}^n$ define $\varphi_\zeta(x)=M\frac{(x\cdot\theta)^2}{2}+x\cdot\zeta.$ 

        Then, there exists $C>0$ and $\tau_0=\tau_0(R_0,V,\lambda)$ such that
        \begin{equation}\label{eq:carlemanest}
          \norm{u}_{X^\half_{-\zeta}}\leq C R_0\norm{e^{\varphi_{\zeta}}(\Delta+\lambda-V)\,(e^{-\varphi_{\zeta}} u)}_{X_{-\zeta}^{-\half}}
        \end{equation}
        for all for $u\in \mathscr{S}(\Rd)$ with $\supp u\subset B_{R_0}$ and all $\zeta$ of the form $\zeta=\tau\theta+i\mathcal{I}$, with $\tau > \tau_0$, $\theta\in\mathbb{S}^{d-1}$ and  $\mathcal{I}\in\Rd$ such that $\mathcal{I}\cdot \theta=0$ and $|\mathcal{I}|\lesssim \tau$.
  \end{prop}

Take now $R_0$ such that $\D\subset B_{R_0}$. We can check that the spaces $X_\zeta^\half$ and $H^1(\Rd)$ are equal as sets, and that, for every $u\in H^1(\Rd)$ such that $\supp u\subset \overline{\D}$, we have that $e^{\varphi_{\zeta}}(\Delta+\lambda-V)\,(e^{-\varphi_{\zeta}} u)$ is in $X_{-\zeta}^{-\half}$. Therefore, by density, the estimate \eqref{eq:carlemanest} also holds for every $u\in H^1(\Rd)$ such that $\supp u\subset \overline{\D}$. This will be useful to prove the following proposition:
\begin{prop}\label{prop:ucp1}
    Consider $d\geq 3$. If $u\in H^1_{\textrm{loc}}(\mathbb{R}^n)$ is a solution of \[ (\Delta+\lambda-V)\,u=0\text{ in } \Rd \] that satisfies the SRC \eqref{eq:src}, then $u$ has to be identically zero.
\end{prop}
\begin{proof}
  Let $R>0$ and call $B=\{x\in\Rd\,:\,|x|<R\}$. On the one hand, the restriction of u to $\Rd\setminus\supp V$ solves $(\Delta+\lambda)\,u=0$. By Theorem 11.1.1 in \cite{hormanderlinear1963} this restriction is smooth, and we have that
  \begin{equation}\label{eq:srcuniq}
    \int_{\partial B} |\partial_\nu u-i\lambda^\half u|^2\,\dd S = \int_{\partial B}\left(|\partial_\nu u|^2+\lambda |u|^2+i\lambda^\half (\partial_\nu u\,\overline{u}-u\,\overline{\partial_\nu u})\right)\dd S,
  \end{equation}
  where $\partial_\nu=\nu\cdot\nabla$ denotes the normal derivative with respect to the vector $\nu=x/|x|$. Using Green's identity in $B\setminus\overline{\D}$ we obtain that
  \begin{equation}\label{eq:srcuniq1}
    \int_{\partial B} \left(\partial_\nu u\,\overline{u}-u\,\overline{\partial_\nu u}\right)\dd S=
    -\int_{\partial \D} \left(\partial_\nu u\,\overline{u}-u\,\overline{\partial_\nu u}\right)\dd S
  \end{equation}
  Now, if we take the limit $R\to\infty$ in \eqref{eq:srcuniq}, the LHS vanishes by the SRC \eqref{eq:src}. This along with \eqref{eq:srcuniq1} yields
  \begin{equation}\label{eq:srcuniq2} \lim_{R\to\infty}\int_{\partial B}\left(|\partial_\nu u|^2+\lambda |u|^2\right) \dd S=i\lambda^\half\int_{\partial \D} \left(\partial_\nu u\,\overline{u}-u\,\overline{\partial_\nu u}\right)\dd S. \end{equation}
  On the other hand, Green's identity in $\D$ gives us that
  \begin{align*}\int_{\partial \D} \left(\partial_\nu u\,\overline{u}-u\,\overline{\partial_\nu u}\right)\dd S&=\int_{\D} \helm u\,\overline{u}-\int_{\D}u\,\overline{\helm u}\\
&= 2i\textrm{Im}\int_{\D}\helm u \overline u =0,\end{align*}
  since $(\Delta+\lambda-V)\,u=0$. Here, $\textrm{Im}$ denotes the imaginary part. This along with \eqref{eq:srcuniq2} implies, in particular, that
  \begin{equation}\label{eq:rellich} \lim_{R\to\infty}\int_{\partial B}\lambda |u|^2=0. \end{equation}
  Now, since $u$ solves $\left(\Delta+\lambda\right)u=0$ in $\Rd\setminus \Omega$ and satisfies the decay condition \eqref{eq:rellich}, Rellich's lemma \cite{Rellich1943} implies that $\supp u\subset \overline{\D}$. This also means that $u\in H^1(\Rd)$. We can therefore take for instace $\zeta = \tau e_n$ with $\tau>\tau_0$ as in Proposition \ref{lem:carlemanest}, and apply inequality \eqref{eq:carlemanest} to $v=e^{\varphi_\zeta}u$, which belongs to $H^1(\Rd)$ and is supported in $\overline{\D}$:
  \begin{equation*}
    \norm{e^{\varphi_{\zeta}}u}_{X^\half_{-\zeta}}\leq C R_0\norm{e^{\varphi_{\zeta}}(\Delta+\lambda-V)\, u}_{X_{-\zeta}^{-\half}},
  \end{equation*}
  where $R_0$ is such that $\D\subset B_{R_0}$. Finally, since $(\Delta+\lambda-V)\,u=0$, we can conclude that $u=0$.
\end{proof}
Another result in the literature for unique continuation of elliptic operators is the following, whose proof can be found in \cite{harrachmonotonicity2019}. We say that an open subset $O\subset \overline \Omega$ is connected to $\Sigma$ if $O$ is connected and $\Sigma\cap O \neq \emptyset$.
\begin{prop}{\cite{harrachmonotonicity2019}}\label{prop:ucp2}
  Let $\Omega$ be a Lipschitz domain, $\Sigma\subset\partial\Omega$ be relatively open and nonempty. Let $C\subset\Omega$ be such that $\Omega\setminus C$ is connected to $\Sigma$. Let also $\lambda>0$. If $u\in H^1(\Omega)$ satisfies
  \begin{equation*}
		\begin{cases}
				\left(\Delta +\lambda\right)u = 0 & \textrm{in }\Omega\setminus C,\\
				u=\partial_\nu u=0 & \textrm{in } \Sigma.
			\end{cases}
	\end{equation*}
  then $u=0$ in $\Omega\setminus C$.
\end{prop}
As an easy consequence of the previous propositions, we have the following:
\begin{cor}\label{prop:ucp}
    Let $\Omega$ be a Lipschitz open domain such that $\supp V\subset \Omega$, and let $\Sigma\subset \partial\Omega$ be relatively open and nonempty. If $u\in H^1(\Omega)$ satisfies
    \begin{equation*}
		\begin{cases}
				\left(\Delta +\lambda-V\right)u = 0 & \textrm{in }\Omega,\\
				u=\partial_\nu u=0 & \textrm{in } \Sigma.
			\end{cases}
	\end{equation*}
  then $u=0$ in $\Omega$.
\end{cor}
\begin{proof}
    Observe that $u$ is a solution to
    \begin{equation*}
		\begin{cases}
				\left(\Delta +\lambda\right)u = 0 & \textrm{in }\Omega\setminus\supp V,\\
				u=\partial_\nu u=0 & \textrm{in } \Sigma.
			\end{cases}
	\end{equation*}
   Since $\supp V$ is closed, $\Omega\setminus \supp V$ is connected to $\Sigma$.
    Therefore, by Proposition \ref{prop:ucp2} u vanishes in $\Omega\setminus\supp V$, and it can be extended by zero to the rest of $\R^n$. If we call $\tilde u$ this extension, we have that $\tilde u\in H^1(\R^n)$ and $\supp \tilde u\subset \supp V$. Finally, if we apply Proposition \ref{prop:ucp1} to $\tilde u$, we obtain $\tilde u\equiv 0$ and the result is proven.
\end{proof}

\subsection{The Scattering Problem}\label{sect:direct}
We will briefly give the definition of the suitable spaces that were used in \cite{Caizares2024} to solve the scattering problem \eqref{eq:scatt}. These were defined in \cite{caroscattering2020} to solve this same problem with potentials of the form $V=V^0+\alpha\dd\sigma$.
\begin{definition}
  Let $Y_\lambda$ be the space of tempered distributions $f\in\mathscr{S}'(\Rd)$ for which the following norm is finite:
  \[\norm{f}^2_{Y_\lambda}\defeq \norm{m_\lambda^{-1/2}\widehat{P_{<I}f}}^2_{L^2}+\sum_{k\in I}\lambda^{-1/2}\norm{P_kf}^2_B+\sum_{k>k_\lambda+1}\norm{m_\lambda^{-1/2}\widehat{P_kf}}^2_{L^2},\]
  where $m_\lambda(\xi)=|\lambda-|\xi|^2|,$ and the norm $\norm{\centerdot}_B$ is defined as
        \[ \norm{f}_{B}=\sum_{j\in\mathbb{N}_0}\,2^{j/2}\norm{f}_{L^2(D_j)},\]
        with $D_j=\{x\in\Rd\;:\;2^{j-1}<|x|\leq2^j\}$, $D_0=\{x\in\Rd:|x|\leq 1\}$.
\end{definition}
\begin{definition}
  Let $Z_\lambda$ be the space of tempered distributions $f\in\mathscr{S}'(\Rd)$ for which the following norm is finite:
  \[\norm{f}^2_{Z_\lambda}\defeq \norm{m_\lambda^{-1/2}\widehat{P_{<I}f}}^2_{L^2}+\sum_{k\in I}\lambda^{n(\frac{1}{q_n'}-\frac{1}{p_n'})}\norm{P_kf}^2_{L^{q_n'}}+\sum_{k>k_\lambda+1}\norm{m_\lambda^{-1/2}\widehat{P_kf}}^2_{L^2},\]
  where $m_\lambda(\xi)=|\lambda-|\xi|^2|$.
\end{definition}
\begin{definition}
  The space $X_\lambda$ is defined as the sum \[X_\lambda\defeq Y_\lambda+Z_\lambda=\{f=g+h\;:\;g\in Y_\lambda,\;h\in Z_\lambda\},\] equipped with the usual norm
  %\[\norm{f}_{X_\lambda}=\inf_{\substack{g\in Y_\lambda\\h\in Z_\lambda\\g+h=f}}\{\norm{g}_{Y_\lambda}+\norm{h}_{Z_\lambda}\}\]
  \[\norm{f}_{X_\lambda}=\inf_{g+h=f}\{\norm{g}_{Y_\lambda}+\norm{h}_{Z_\lambda}\}.\]
\end{definition}
\begin{definition}
  The Banach space $(X_\lambda^*,\norm{\centerdot}_{X_\lambda^*})$ is defined as the dual space of \newline$(X_\lambda,\norm{\centerdot}_{X_\lambda})$.
\end{definition}
Actually, the space $X_\lambda^*$ is isomorphic to the space of $u\in\mathscr{S}'(\Rd)$ for which the following norm is finite:
\begin{align*}
  \norm{u}^2_{X^*_\lambda}\defeq& \norm{m_\lambda^{1/2}\widehat{P_{<I}f}}^2_{L^2}+\sum_{k\in I}\left(\lambda^{1/2}\norm{P_kf}^2_{B^*}+\lambda^{n(\frac{1}{q_n}-\frac{1}{p_n})}\norm{P_kf}^2_{L^{q_n}}\right)+\\&\sum_{k>k_\lambda+1}\norm{m_\lambda^{1/2}\widehat{P_kf}}^2_{L^2},
\end{align*}
where the norm $\norm{\centerdot}_{B^*}$ is defined by
\[\norm{u}_{B^*}=\sup_{j\in\mathbb{N}_0}\,\left(2^{-j/2}\norm{f}_{L^2(D_j)}\right)\]
Furthermore, it might interesting to note that $\mathscr{S}(\Rd)$ is dense in both $X_\lambda$ and $X_\lambda^*$. It will later be useful to see the elements of $X_\lambda^*$ as elements of $H^1_{\textrm{loc}}(\R^n)$, and in fact we have the following bound:
\begin{prop}[\cite{Caizares2024}]\label{prop:restriction}
  For any  bounded open domain $\Omega\subset\Rd$, the restriction map
  \begin{align*}
    r_\Omega\,:\,X_\lambda^*&\longrightarrow H^1(\Omega)\\
    u&\longmapsto\restr{u}{\Omega}
  \end{align*}
  is a bounded operator.
\end{prop}
The following theorem gives the existence of a unique solution to the scattering problem:
\begin{theorem}[\cite{Caizares2024}]\label{thm:direct}
  Suppose $n\geq3$ and $V$ is of the form \eqref{V}. Then, there exists $\lambda^V=\lambda^V(V,n)$ such that, for every $\lambda\geq \lambda^V$, there is a unique solution $u\in X_\lambda^*$ to the problem 
  \[
  \begin{cases}
    (\Delta+\lambda-V)u=f\textnormal{ in }\Rd,\\
    u\textnormal{ satisfies SRC},
  \end{cases}  
  \]
  for every $f\in X_\lambda$. Moreover, the mapping $f\mapsto u$ is bounded from $X_\lambda$ to $X_\lambda^*$.
\end{theorem}

\section{Layer potentials}\label{sect:layer}
Layer potentials are a classical tool in the study of elliptic boundary value problems. They allow us to turn these problems into integral equations with kernels associated to the fundamental solution. According to \cite{Steinbach}, the idea was introduced by Gauss in 1839, and then developed by Neumann in the 1870s and 1880s for convex domains. It was then studied by Poincaré in the case of smooth domains. Several generalizations have been made over the years. We can cite for instance the recent works of David and Semmes \cite{MR1113517} and of Tolsa \cite{tolsajump2020}, where they study general layer potentials defined by weakly singular kernels on rectifiable sets. However, we restrict our attention to layer potentials defined by the fundamental solution for the operator $\Delta+\lambda-V$ in $\C^2$ domains. Our main references here are the book by Folland \cite{follandintroduction1995}, who studies the Laplacian on $\mathcal \C^2$ domains in $\R^n$, and that of Colton and Kress \cite{coltoninverse2013}, who focus on the Helmholtz equation in $\mathcal{C}^2$ domains in $\R^3$, as well as the book \cite{coltonintegral} by the same authors. We also point to the books \cite{Courant2008} and \cite{Kellogg1967} for classical references.
Along this section, let $V$ be a potential of the type \eqref{V} and fix $\lambda>\lambda^V$, where $\lambda^V>0$ is the lower bound required by Theorem \ref{thm:direct} for the solvability of the forward problem. From now on, fix also a bounded open domain $\Omega$ of class $\C^2$ such that $\supp V\subset\Omega$. 

Remember that we want to study the following general Neumann boundary: find $u\in H^1(\Omega)$ such that
\begin{equation}\label{eq:neumanngen}
    \begin{cases}
            \left(\Delta +\lambda- V\right)u= f_1 & \textrm{in }\Omega,\\
            \partial_{\nu}u=f_2 & \textrm{on }\partial \Omega,
        \end{cases}
\end{equation}
for $f_1\in H^{-1}_0(\Omega)\cap L^2(\Omega\setminus\supp V)$ and $f_2\in L^2(\partial\Omega)$. In this section, we are going to characterize the pairs $(f_1,f_2)$ that admit a solution in terms of the eigenspace associated to $\lambda$ and, in particular, we will show that one can solve the problem above whenever $\lambda$ is not a NEV for $-\Delta+V$ in $\Omega$.

Let's start by doing a transformation on the equation. Let $w\in X_\lambda^*$ be the unique solution to 
\begin{equation}\label{eq:w}
  \begin{cases}
      \left(\Delta +\lambda- V\right)w = f_1 & \textrm{in }\Rd,\\
      w\textrm{ satisfying SRC},
    \end{cases}
  \end{equation}
  which exists by Theorem \ref{thm:direct}, since the extension by $0$ of $f_1$ to $\R^n$ lives in $X_\lambda$. Remember that the spaces $X_\lambda$ and $X_\lambda^*$ were defined in Section \ref{sect:direct}. Denote now $g=f_2-\restr{\partial_\nu w}{\partial\Omega}$, and note that, since $w$ is a solution to $(\Delta+\lambda)w=f_1$ away from $\supp V$, then $w\in H^2(\R^n\setminus\supp V)$ and therefore $\restr{\partial_\nu w}{\partial\Omega}\in L^2(\Omega)$. Thus solving \eqref{eq:neumanngen} is equivalent to solving the problem 
  \begin{equation}\label{eq:neumann2}
    \begin{cases}
            \left(\Delta +\lambda- V\right)v= 0 & \textrm{in }\Omega,\\
            \partial_{\nu}v=g & \textrm{on }\partial \Omega.
        \end{cases}
\end{equation}
with $g\in L^2(\Omega)$ by setting $u=w+v$. 

We will study this problem using the method of layer potentials. Note also that $\restr{w}{\Omega}\in H^1(\Omega)$ by Proposition \ref{prop:restriction}, so we will only need to prove that $v$ belongs to $H^1(\Omega)$ to conclude that $u$ belongs too.
In this case it will be useful to think of the fundamental solution for the operator $\Delta+\lambda$ in $\Rd$ as a Hankel function. In fact, $\Phi_\lambda$ as defined by \eqref{eq:homo} will take the form
\begin{equation}
    \Phi_\lambda(x)=\frac{i}{4}\left(\frac{\lambda^\half}{2\pi|x|}\right)^{n/2-1}H_{n/2-1}^{(1)}\left(\lambda^\half |x|\right),
\end{equation}
where $H_\nu^{(1)}$ denotes the Hankel function of the first kind (or Bessel function of the third kind). If we define $u_{in}(x,y)=\Phi_\lambda(x-y)$ and recall the limiting properties of the Hankel functions, it's relatively easy to check that
\begin{align}\label{eq:kernel}
    \begin{split}
        u_{in}(x,y)&= F(x,y)|x-y|^{2-n},\\
        \partial_{\nu_x}u_{in}(x,y)&=G(x,y)|x-y|^{2-n},\\
        \partial_{\nu_y}u_{in}(x,y)&=\partial_{\nu_x}u_{in}(y,x)=G(y,x)|x-y|^{2-n},
    \end{split}
\end{align}
with $F$ and $G$ being bounded functions on $\partial \Omega\times\partial \Omega$. Then, $u_{in}$, $\partial_{\nu_x}u_{in}$ and $\partial_{\nu_y}u_{in}$ are, by definition, weakly singular kernels of order $n-2$ on $\partial \Omega\times\partial \Omega$. The following lemma is a combination of those in \cite{follandintroduction1995}, Chapter 3B, and is an important piece to solve the problem \eqref{eq:neumann2}.
\begin{lemma}[\cite{follandintroduction1995}]\label{lem:kernel}
    If we denote by $T$ the integral operator over $\partial \Omega$ defined by a weakly singular kernel $K$ of order $\alpha$ on $\partial \Omega\times \partial \Omega$, with $0<\alpha<n-1$, as
    \[ \left(Tf\right)(x)=\int_{\partial \Omega}T(x,y)\,f(y)\dd S(y), \] then the following statements hold:
    \begin{enumerate}
        \item $T$ is compact on $L^2(\partial \Omega)$,
        \item $T$ transforms bounded functions into continuous functions, and
        \item if $f\in L^2(\partial \Omega)$ and $f+Tf\in \mathcal{C}(\partial \Omega)$, then $f\in\mathcal{C}(\partial \Omega)$.
    \end{enumerate}
\end{lemma}
Take now $u_{to}=u_{in}+u_{sc}$, where $u_{sc}$ is the scattering solution of \eqref{eq:scatt}, and define for $f$ continuous on $\partial \Omega$ and $x\in \Rd\setminus \partial \Omega$ the \textbf{single layer potential} with moment $f$ as
\[ \left(\mathcal{S}f\right)(x)=\int_{\partial \Omega}u_{to}(x,y)\,f(y)\dd S(y), \]
and the \textbf{double layer potential} as 
\[ \left(\mathcal{D}f\right)(x)=\int_{\partial \Omega}\partial_{\nu_y}u_{to}(x,y)\,f(y)\dd S(y),. \]
Define further the operator $\mathcal{N}$, which is the adjoint of $\mathcal{D}$ over $\partial \Omega$, as
\[ \left(\mathcal{N}f\right)(x)=\int_{\partial \Omega}\partial_{\nu_x}u_{to}(x,y)\,f(y)\dd S(y),\quad x\in \partial \Omega, \]
which must be understood as an improper integral.

Now, we have the following lemmas, which are modifications of classical results that can be found in \cite{coltoninverse2013}, \cite{follandintroduction1995} and \cite{coltonintegral} for the Helmholtz and Laplace equations. We denote by $u_+$ the trace on $\partial \Omega$ of $\restr{u}{\Rd\setminus \overline{\Omega}}$, and by $u_-$ the trace on $\partial \Omega$ of $\restr{u}{\Omega}$. As well, we denote by $\partial_\nu u_+$ and $\partial_\nu u_-$ the normal derivative of those, always with respect to the outward-pointing normal vector of $\partial \Omega$ (as seen from inside $\Omega$).
\begin{lemma}\label{lem:single}
    Consider $n\geq3$. Lef $f\in\mathcal{C}(\partial \Omega)$. Then, the single layer potential $u=\mathcal{S}f$  is continuous throughout $\Rd$, and we have the limiting values
    \begin{equation}\label{eq:bdsingle}
      \partial_\nu u_{\pm}(x)=\left(\mathcal{N}f\right)(x)\mp \frac{1}{2}f(x),\quad x\in\partial \Omega,
  \end{equation}
  where the integral exists as an improper integral.
  Consequently, we have the jump relation $\partial_\nu u_--\partial_\nu u_+ = f$ on $\partial \Omega$. Furthermore, $u$ is a solution in $H^1_{\textrm{loc}}(\Rd)$ to $\helm u = 0$ in $\Rd\setminus \partial \Omega$ and fullfils SRC \eqref{eq:src}. 
\end{lemma}
\begin{proof}
    First, note that the single layer potential for the homogeneuos Helmholtz equation \[ v(x)=\int_{\partial \Omega}u_{in}(x,y)f(y)\dd S(y) \] can be extended to the boundary, is a solution in $H^1_{\textrm{loc}}(\Rd)$ to $\left(\Delta+\lambda\right) v= 0$ in $\Rd\setminus \partial \Omega$, fullfils SRC \eqref{eq:src} and has boundary values     
    \begin{equation*}
        \partial_\nu v_{\pm}(x)=\int_{\partial \Omega}\partial_{\nu_x}u_{in}(x,y)\,f(y)\dd S(y)\mp \frac{1}{2}f(x),\quad x\in\partial \Omega,
    \end{equation*} which is a classical result, see for example \cite{coltoninverse2013}.
    Now, define \[ w(x)=\int_{\partial \Omega}u_{sc}(x,y)f(y)\dd S(y). \] To see that $w$ is in $H^1_{\textrm{loc}}(\Rd)$, take $K\in\Rd$ compact, then
    \[ \norm{w}_{H^1(K)}\lesssim \sup_{y\in\partial\Omega}\norm{u_{sc}({\centerdot\,,y})}_{H^1(K)}\leq\sup_{y\in\partial\Omega} \norm{u_{sc}({\centerdot\,,y})}_{X_\lambda^*} \lesssim \sup_{y\in\partial\Omega}\norm{Vu_{in}({\centerdot\,,y})}_{X_\lambda}. \]
    In the second inequality, we have used Proposition \ref{prop:restriction}, while the third one is the result of Theorem \ref{thm:direct}. For the first one, observe that
    \[\begin{aligned}
      \norm{w}^2_{L^2(K)}&=\int_{K}\left|\int_{\partial \Omega}u_{sc}(x,y)f(y)\dd S(y)\right|^2\dd x\\&\leq \int_{K}\left(\int_{\partial \Omega}|u_{sc}(x,y)|^2\,|f(y)|^2\dd S(y)\right)\dd x\\
      &=\int_{\partial \Omega}\left(\int_{K}|u_{sc}(x,y)|^2\dd x\right)|f(y)|^2\dd S(y)\\&= \int_{\partial \Omega}\norm{u_{sc}(\centerdot\,,y)}_{L^2(K)}^2\,|f(y)|^2\dd S(y)\\
      &\leq \norm{f}_{L^2(\partial\Omega)}^2\sup_{y\in\partial\Omega}\norm{u_{sc}({\centerdot\,,y})}^2_{L^2(K)},
    \end{aligned}\]
    where we have used Fubini's Theorem and the Cauchy-Schwarz inequality. Then, to get the inequality in the $H^1(K)$ norm, we have to differentiate under the integral sign.
    If we take a smooth cut-off function $\eta$ such that $\eta\equiv 1$ in $\supp V$ and $\eta\equiv 0$ in $\partial\Omega$, we have that $V(x)u_{in}(x,y)=V(x)\eta(x)u_{in}(x,y)$ and
    \[ \norm{Vu_{in}({\centerdot\,,y})}_{X_\lambda}\lesssim \norm{\eta u_{in}(\centerdot\,,y)}_{X_\lambda^*}\lesssim 1, \]
    where we have used that multiplication by $V$ is bounded from $X_\lambda^*$ to $X_\lambda$, as recalled in Section \ref{sect:direct}, and that $u_{in}(\centerdot\,,y)$ is smooth away from $y$ by Theorem 11.1.1 in \cite{hormanderlinear1963}, since it solves $(\Delta+\lambda)u_{in}(\centerdot\,,y)=0$. This proves that $w$ is in $H^1_{\textrm{loc}}(\Rd)$.

     Moreover, since $u_{sc}$ solves the problem \eqref{eq:scatt}, it is easy to check that $u = v + w$ solves $\helm u = 0$ in $\Rd\setminus \partial \Omega$ and fullfils the SRC \eqref{eq:src}. Also, since for any $y\in\partial \Omega$, $u_{sc}(\centerdot\,,y)$ solves $\left(\Delta+\lambda\right) u_{sc}(\centerdot\,,y)= 0$ in $\Rd\setminus\supp V$, it is smooth in this set by Theorem 11.1.1 in \cite{hormanderlinear1963}, and in particular it is smooth near $\partial \Omega$. Therefore, the limiting values of $w$ on the boundary are just
    \[ \partial_\nu w_{\pm}(x)=\int_{\partial \Omega}\partial_{\nu_x}u_{sc}(x,y)\,f(y)\dd S(y),\quad x\in\partial \Omega, \] and therefore \eqref{eq:bdsingle} is fullfiled. 
\end{proof}
\begin{lemma}\label{lem:double}
    Consider $n\geq3$. Lef $f\in\mathcal{C}(\partial \Omega)$. Then, the double layer potential $u=\mathcal{D}f$  can be extended continuously to $\partial \Omega$, and we have the limiting values
    \begin{equation}\label{eq:bddouble}
        u_{\pm}(x)=\left(\mathcal{D}f\right)(x)\pm \frac{1}{2}f(x),\quad x\in\partial \Omega,
    \end{equation}
    where the integral exists as an improper integral. Consequently, we have the jump relation $u_+-u_- = f$ on $\partial \Omega$. Furthermore, $u$ a solution in $H^1_{\textrm{loc}}(\Rd)$ to $\helm u = 0$ in $\Rd\setminus \partial \Omega$, it fullfils SRC \eqref{eq:src} and {${\partial_\nu u_--\partial_\nu u_+=0}$} on $\partial\Omega$.
\end{lemma}
\begin{proof} 
    The proof goes exactly as the proof of Lemma \ref{lem:single} above, we just need to make a comment on how to prove the last statement. Indeed, the fact that the double layer potential for the homogeneous Helmholtz equation \[ v(x)=\int_{\partial \Omega}\partial_{\nu_y} u_{in}(x,y)f(y)\dd S(y) \] fulfills that $\partial_\nu v_--\partial_\nu v_+=0$ on $\partial\Omega$ is classical (see for example \cite{coltoninverse2013}). Meanwhile, the function \[ w(x)=\int_{\partial \Omega}\partial_{\nu_y}u_{sc}(x,y)f(y)\dd S(y) \] is smooth away from $\supp V$, since $u_{sc}(\centerdot\,,y)$ is smooth away from $\supp V$ as well. 
\end{proof}
Now, we can try to solve the problem \eqref{eq:neumann2} by means of Fredholm alternative. This is done in the following lemma.
\begin{lemma}\label{lem:np1}
    Let $g\in L^2(\partial\Omega)$. There exists a solution $v\in H^1(\Omega)$ to the Neumann problem
    %\begin{equation}\label{eq:neu}
    %  \begin{cases}
    %    \helm v = 0 & \textup{ in } \Omega,\\
    %    \partial_\nu v = g &  \textup{ on } \partial\Omega,
    %  \end{cases}
    %\end{equation}
    \eqref{eq:neumann2} if and only if
    \[ \int_{\partial \Omega} g\, u =0\]
    for all $u\in H^1(\Omega)$ solution to the homogeneous Neumann problem
    \begin{equation}\label{eq:neuhom}
      \begin{cases}
        \helm u = 0 & \textup{ in } \Omega,\\
        \partial_\nu u = 0 &  \textup{ on } \partial\Omega,
      \end{cases}
    \end{equation}
  \end{lemma}
  \begin{proof}
    The necessity of the condition can be seen using integration by parts. Indeed, if there exists such a solution $v$, then
    \begin{align*}
      0=\int_{\Omega}\helm v\, u = \int_{\Omega}v\helm u + \int_{\partial\Omega}\partial_\nu v\,u-\int_{\partial\Omega}v\,\partial_\nu u 
      =\int_{\partial\Omega}g\,u.
    \end{align*}
    Now, to prove that the condition is sufficient, observe that, by Lemma \ref{lem:single}, if we can find $\phi\in L^2(\partial \Omega)$ solving the integral equation $\mathcal{N}\phi+\frac{1}{2}\phi=g$, then $v=\mathcal{S}\phi$ is a solution to \eqref{eq:neumann2}. 
    
    Note that both $\mathcal{N}$ and $\mathcal{D}$ are compact operators on $L^2(\partial\Omega)$, and it holds that $\mathcal{D}=\mathcal{N}^*$. Therefore, by Fredholm alternative (see e.g. theorem 1.29 in \cite{coltonintegral}), we have that $\textrm{ran}(\mathcal{N}+\frac{1}{2}I)=\textrm{ker}(\mathcal{D}+\frac{1}{2}I)^\perp$.

    Set now $\psi\in\ker(\mathcal{D}+\frac{1}{2}I)$ and let $u=\mathcal{D}\psi$. Then, $u$ is a solution to the problem
    \begin{equation*}
    \begin{cases}
    \left(\Delta+\lambda\right) u = 0 & \textrm{in }\Rd\setminus \overline{\Omega},\\
    u_+ = 0 & \textrm{on }\partial \Omega,\\ 
    u\textrm{ satisfying SRC},
    \end{cases}
  \end{equation*}
    Therefore, $u=0$ in $\Rd\setminus \Omega$, by uniqueness of the exterior Dirichlet problem \cite{coltonintegral}. Hence, $\partial_\nu u=0$, as the normal derivative of the double layer potential is continous through the boundary. This means that $u$ is as well a solution to the problem \eqref{eq:neuhom},
  and $\psi=u_+-u_-=-u_-$, which means that \[\ker(\mathcal{D}+\frac{1}{2}I)=\{\restr{u}{\partial\Omega}:u\textup{ solution to the homogeneous Neumann problem \eqref{eq:neuhom}}\}, \]
  and therefore the condition is sufficient.
  \end{proof}
  \begin{cor}\label{cor:finitenev}
    The eigenspace of any Neumann eigenvalue $\lambda$ for $-\Delta+V$ in $\Omega$ is finite-dimensional.
  \end{cor}
  \begin{proof}
    Let $u_1$ and $u_2$ be two solutions to the homogeneous Neumann problem such that $u_1=u_2$ in $\partial\Omega$. Then $v=u_1-u_2$ solves
    \begin{equation*}
		\begin{cases}
				\left(\Delta +\lambda-V\right)v = 0 & \textrm{in }\Omega,\\
				v=\partial_\nu v=0 & \textrm{in } \partial\Omega.
			\end{cases}
	\end{equation*}
     which means that $v=0$ in $\Omega$ by the unique continuation property, Proposition \ref{prop:ucp}. Therefore, two distinct NEVs for $\lambda$ will also have distinct boundary data, and thus the eigenspace of $\lambda$ coincides with the kernel of the compact operator $\mathcal{D}+\frac{1}{2}I$ by the previous proof, and is therefore finite-dimensional.
  \end{proof}
  \begin{cor}\label{cor:np2}
    Let $f_1\in H^{-1}_0(\Omega)\cap L^2(\Omega\setminus\supp V)$ and $f_2\in L^2(\partial\Omega)$. There exists a solution $u\in H^1(\Omega)$ to the Neumann problem \eqref{eq:neumanngen}
  %\begin{equation}\label{eq:neu2}
  %  \begin{cases}
  %    \helm v = f_1 & \textup{ in } \Omega,\\
  %    \partial_\nu v = f_2 &  \textup{ on } \partial\Omega,
  %  \end{cases}
  %\end{equation}
  if and only if
  \[ \int_{\partial \Omega} f_2\, v = \int_{\Omega} f_1\, v\]
  for all $v\in H^1(\Omega)$ solution to the homogeneous Neumann problem \eqref{eq:neuhom}.
  \end{cor}
  \begin{proof}
  Let $w\in X_\lambda^*$ be the unique solution to the equation 
  \begin{equation*}
    \begin{cases}
        \left(\Delta +\lambda- V\right)w = f_1 & \textrm{in }\Rd,\\
        w\textrm{ satisfying SRC},
      \end{cases}
    \end{equation*}
    denote $g=f_2-\restr{\partial_\nu w}{\partial\Omega}$. There exists a solution $u=\tilde u-w$ to the problem \eqref{eq:neumanngen} if and only if $\tilde u$ solves \eqref{eq:neumann2}. Using Lemma \ref{lem:np1}, and integrating by parts in $\Omega$, we get to the sufficient and necessary condition
    \begin{align*}
      0=&\int_{\partial\Omega}\left(f_2-\partial_\nu w\right) v\\
      =&\int_{\partial\Omega}f_2 \,v-\int_{\Omega}\helm w\,v+\int_{\Omega} w\,\helm v+\int_{\partial\Omega}w\,\partial_\nu v\\
      =&\int_{\partial\Omega}f_2 \,v-\int_{\Omega}f_1\,v
    \end{align*}
    for all $v\in H^1(\Omega)$ solution to the homogeneous Neumann problem in $\Omega$.
\end{proof}
Finally, we can conclude the following result:
\begin{cor}\label{cor:nev}
    There exists a unique solution to the problem \eqref{eq:neumanngen} for every $f_1\in H^{-1}_0(\Omega)\cap L^2(\Omega\setminus\supp V)$ and $f_2\in L^2(\partial\Omega)$ if and only if $\lambda>\lambda^V$ is not a Neumann eigenvalue for $-\Delta+V$ in $\Omega$.
\end{cor}

\section{Perturbation of domains}\label{sect:calculus}
\subsection{The approach of $\C^m$ embeddings}

Let $V$ be a potential as in \eqref{V}, $\lambda\geq 0$ and take $\Omega$ and be a $\C^m$ domain such that and $\supp V\subset\Omega$. From now on, we will denote $\supp V$ by $\Vcal$. Define then the operator
\begin{align*}L_\Omega : H^1(\Omega)\cap H^m(\Oprime)&\to H^{-1}_0(\Omega)\cap H^{m-2}(\Oprime)\\
u & \mapsto \left(\Delta+\lambda-V\right) u.\end{align*}
We will deform the domain $\Omega$ by a \textbf{$\mathcal{C}^m$ embedding} $h:\Omega\to\Rd$. This is, $h$ is of class $\C^m$, it is a diffeomorphism onto its image and its inverse $h^{-1}$ is also of class $\C^m$. Denote the space of such embeddings as $\textrm{Diff}^m(\Omega)$. We would like to study the natural operator that acts on functions on the deformed domain
\[L_{h(\Omega)} : H^1(h(\Omega))\cap H^m(h(\Oprime))\to H^{-1}_0(h(\Omega))\cap H^{m-2}(h(\Oprime)).\] However, the fact that the function spaces change with $h$ poses a difficulty. To solve this, we will consider the pullback $h^*$, defined by \[h^*u(x)=u\, (h(x)),\] and the pushforward $h_*$, defined by \[h_*u(x)=u \,(h^{-1}(x)).\] Then, the \textit{Lagrangian form} of the operator $L_{h(\Omega)}$ is defined as \[h^* L_{h(\Omega)} h_*: H^1(\Omega)\cap H^m(\Oprime)\to H^{-1}_0(\Omega)\cap H^{m-2}(\Oprime).\] 

Take now a curve of class $\C^1$ of such embeddings $t\mapsto h_t$. In \cite{Henry2005}, Henry introduces the anti-convective derivative, defined for $\C^1$ functions $\varphi:\R\times\R^n\to \R$ as
\begin{equation}\label{eq:anti-convective}
  D_t\varphi(t,x)=\partial_t \varphi(t,x)-U(t,x)\cdot\nabla \varphi(t,x),
\end{equation}
where \[U(t,x)=\left((Jh_t)^{-1}\,\partial_t h_t\right)(x)\] and $Jh$ denotes the Jacobian matrix of $h$ with respect to $x$. This object will allow us to differentiate operators and boundary conditions in the Lagrangian form. It satisfies the following rule, which can be proved as consequence of the chain rule:
\begin{lemma}\label{lem:convective}
  Suppose that $\psi:\R\times\R^n\to\R$ and $h:\R\times\R^n\to\R^n$ are of class $\mathcal{C}^1$, and call $h_t=h(t,\centerdot)$.
  
  Then, for those $(t,x)\in\R\times\R^n$ such that ${\textnormal{det }(Jh)(t,x)\neq0}$ we have
  \[D_t(h_t^*\,\psi(t,\centerdot))(x)=(h_t^*\,\partial_t\psi(t,\centerdot))(x),\] 
  where $h_t^*$ is the pullback by $h_t$.
\end{lemma}
Although stated pointwise, this lemma and the results below will be later interpreted as an equality between $\C^1$ curves with values on Sobolev spaces. This interpretation is inmediate for curves valued in $\C^1$, and can be extended to Sobolev spaces by density. 

We will now use this to differentiate the differential operator $\Delta+\lambda-V$ with respect to boundary deformation. The reader can be pointed to Chapter 2 of \cite{Henry2005} for an expression of the derivative of a general class of differential operators.

Consider that $t\mapsto \lambda(t)$, $t\mapsto u(t,\centerdot)$ and $t\mapsto h_t$ are curves of class $\C^1$. Suppose that $h_0=i_\Omega$, the trivial embedding of $\Omega$ into $\R^n$. We want to find now
\[\restr{\partial_t\big(h_t^*L_{\Omega_t}h_{t*}u(t,\centerdot)\big)}{t=0},\]
where $\Omega_t=h_t(\Omega)$. Define then $v(t,x)=h_{t*}u(t,x)$ and $\psi(t,x)=L_{\Omega_t}v(t,x)$. Then, by the definition of the anticonvective derivative \eqref{eq:anti-convective} and Lemma \ref{lem:convective}, we have formally that
\begin{align*}
  \partial_t\big(h_t^*\psi(t,\centerdot)\big)=&D_t\left(h_t^*\psi(t,\centerdot)\right) + U(t,\centerdot)\cdot\nabla (h_t^*\psi(t,\centerdot))\\
  =&h_t^*\partial_t\psi(t,\centerdot) + U(t,\centerdot)\cdot\nabla (h_t^*\psi(t,\centerdot))\\
  =&h_t^*\partial_tL_{\Omega_t}v(t,\centerdot)+h_t^*L_{\Omega_t}\,\partial_tv(t,\centerdot) + U(t,\centerdot)\cdot\nabla (h_t^*\psi(t,\centerdot)).
\end{align*}
Note now that, again by Lemma \ref{lem:convective}
\[h_{t*}D_t u = h_{t*}D_t\left(h_t^*v(t,\centerdot)\right)= h_{t*}h_t^*\partial_t v(t,\centerdot)=\partial_t v(t,\centerdot).\]
Therefore, if we observe that $\partial_tL_{\Omega_t}=\lambda'(t)$, and using again the definition of the anticonvective derivative for $u$, we obtain
\begin{equation}\label{eq:derivativeL}
  \begin{aligned}
    \restr{\partial_t\big(h_t^*L_{\Omega_t}h_{t*}u(t,\centerdot)\big)}{t=0}=&L_\Omega(\dot u-\dot h \cdot\nabla u)+\dot\lambda u+\dot h \cdot \nabla (L_\Omega u)\\
    =&L_\Omega \dot u + \dot\lambda u + \big[\dot h \cdot \nabla,L_\Omega\big].
  \end{aligned}
\end{equation}
Here, we have denoted with $\dot u$, $\dot \lambda$ and $\dot h$ the respectives partial derivaties of $u$, $\lambda$ and $h$ with respect to $t$ at $t=0$, and we have used the fact that $Jh_0=Ji_\Omega=1$, so that $U(0,\centerdot)=\dot h$.

Note that this expression may not make sense in general, as the quantity
\[\big[\dot h \cdot \nabla,L_\Omega\big]=\big[\dot h \cdot \nabla,V\big]\] 
poses technical problems due to the low regularity of the potential $V$. However, we will later restrict our attention to curves of deformations of the form
\[h_t(x)=x+t\X(x),\]
for vector fields $\X\in \C^k(\Omega;\R^n)$, $k\geq 1$, such that $\supp \X\cap\supp V = \emptyset$. Note that $h_t\in\textup{Diff}^k(\Omega)$ as long as $Jh_t\neq0$, i.e. for $t$ small enough. 

We also need to be able to differentiate the Neumann boundary conditions. For this purpose, we choose an extension of $\nu_\Omega$ near $\partial\Omega$, where $\nu_\Omega$ is the unit outward-pointing normal vector to $\partial\Omega$, and then define $\nu_{\Omega_t}=\nu_{h_t(\Omega)}$ by the expression
\begin{equation}\label{eq:hnormal}
  h^*\nu_{\Omega_t}(x)=\frac{((Jh_t)^{-1})^\top\nu_\Omega(x)}{|((Jh_t)^{-1})^\top\nu_\Omega(x)|}
\end{equation}
for $x$ near $\partial\Omega$. Here $((Jh)^{-1})^\top$ is the inverse-transpose of the Jacobian matrix of $h$ with respect to $x$. An expression for the derivative of this extension of the normal vector is given in the following lemma from \cite{Henry2005}.
\begin{lemma}
  Let $\Omega$ be a domain of class $\C^2$, and for $h_t\in \textup{Diff}^2(\Omega)$ define $\nu_{\Omega_t}$ as in \eqref{eq:hnormal}. Then, if $t\mapsto h_t$ is a $\C^1$ curve of embeddings of the form $h_t(x)=x+t\X(x)$, with $\X\in\C^2(\Omega;\R^n)$, we can compute the derivative of $\nu_{\Omega_t}$ as
  %\[\frac{\partial\nu_{\Omega_t}}{\partial t}(y)=D_t(h_t^*\nu_{\Omega_t})(x)=-\big(\nabla_{\partial\Omega_t}\sigma+\sigma\frac{\partial\nu_{\Omega_t}}{\partial\nu_{\Omega_t}}\big)(y)\]
  \[\frac{\partial\nu_{\Omega_t}}{\partial t}(y)=D_t(h_t^*\nu_{\Omega_t})(x)=-\nabla_{\partial\Omega_t}\big(\X\cdot\nu_{\Omega_t})(y)\]
  for $y=h_t(x)$ near $\partial\Omega_t$, $x$ near $\partial\Omega$.%, where $\sigma=\X\cdot\nu_{\Omega_t}$ is the normal velocity.
\end{lemma}
%Note that, in the sequel, we will take an extension of the normal vector such that $\partial\nu_{\Omega_t}/\partial\nu_{\Omega_t}=0$ near $\partial\Omega$.
%\begin{theorem}
%  Suppose $b(t,y,\lambda,\mu)$ be of class $\C^1$ in an open set of $\R\times\R^n\times\R^p\times\R^q$, $L,M$ are constant-coefficient differential operators of order $\leq m$ with $Lv(y)\in\R^p$, $Mn(y)\in\R^q$. Assume $\Omega\subset \R^n$ is a $\C^{m+1}$ region, define for $h\in \textup{Diff}^{m+1}(\Omega)$ the operator
%  \[B_{h(\Omega)}(t)v(y)=b\big(t,y,Lv(y),M\nu_{h(\Omega)}(y)\big),\]
%  for $y\in h(\Omega)$ near $\partial h(\Omega)$, where $\nu_{h(\Omega)}$ is the normal vector to $\partial h(\Omega)$ defined as in \eqref{eq:hnormal}. 
%  
%  If $t\mapsto h_t$ be a $\C^1$ curve of embeddings in $\textup{Diff}^m(\Omega)$, then at points in $\Omega$ near $\partial \Omega$ we have
%  \[\begin{aligned}D_t\big(h^*B_{\Omega_t}(t)h_*\big)(u)=&\big(h^*\dot B_{\Omega_t}(t)h_*\big)(u)+\big(h^*B'_{\Omega_t}(t)h_*\big)(u)\cdot D_tu\\
%  &+\big(h^*\frac{\partial B_{\Omega_t}}{\partial\nu}(t)h_*\big)(u)\cdot D_t\big(h^*\nu_{\Omega_t}\big),\end{aligned}\]
%  where
%  \[\dot B_Q(t)v(y)=\frac{\partial b}{\partial t}(t,y,Lv(y)), \]
%  and
%  \[B'_Q(t)v(y)\cdot w(y)=\frac{\partial b}{\partial \lambda}(t,y,Lv(y))Lw(y)\]
%  is the linearization of $v\mapsto B_Q(t)v$, and
%  \[\frac{\partial B_{\Omega_t}}{\partial\nu}(t)(v)\cdot\eta(y)=\frac{\partial b}{\partial\mu}(t)\big(t,y,Lv(y),M\nu_{h(\Omega)}(y)\big)\cdot M\eta(y).\]
%\end{theorem}

This allows us to do a similar computation as before for the Neumann boundary conditions. Indeed, we want to find an expression for
\[\restr{\partial_t\big(h_t^*\nu_{\Omega_t}\cdot\nabla (h_{t*}u)\big)}{t=0}\]
in points of $\partial\Omega$. Consider a curve of embeddings $t\mapsto h_t$ of the form \linebreak $h_t(x)=x+t\X(x)$, with $\X\in\C^2(\Omega;\R^n)$. Define again $v(t,x)=h_{t*}u(t,x)$ and $\phi(t,x)=\nu_{\Omega_t}\cdot\nabla v(t,x)$. Then, by the definition of the anticonvective derivative \eqref{eq:anti-convective} and Lemma \ref{lem:convective}, we have that
\begin{align*}
  \partial_t \psi(t,\centerdot)&=\nu_{\Omega_t}\cdot\nabla\big(\partial_t v(t,\centerdot)\big)+\partial_t\nu_{\Omega_t}\cdot\nabla v(t,\centerdot)\\
  &=\nu_{\Omega_t}\cdot\nabla(h_{t*}D_t u)-\nabla_{\partial\Omega_t}\big(\X\cdot\nu_{\Omega_t}\big)\cdot\nabla v,
\end{align*}
and thus
\begin{align*}
  \partial_t\big(h_t^*\psi(t,\centerdot)\big)=&D_t\left(h_t^*\psi(t,\centerdot)\right) + U(t,\centerdot)\cdot\nabla (h_t^*\psi(t,\centerdot))\\
  =& h_t^*\partial_t\psi(t,\centerdot) + U(t,\centerdot)\cdot\nabla (h_t^*\psi(t,\centerdot))
\end{align*}
means that
\begin{equation}\label{eq:derivativeN}
    \restr{\partial_t\big(h_t^*\nu_{\Omega_t}\cdot\nabla (h_{t*}u)\big)}{t=0}= \partial_\nu\big(\dot u-\dot h \nabla u\big)-\nabla_{\partial\Omega_t}\big(\dot h\cdot\nu\big)\cdot\nabla u+\dot h\cdot\nabla\big(\partial\nu u\big)
\end{equation}
Here we have again used the fact that $Jh_0=Ji_\Omega=1$, so that $U(0,\centerdot)=\dot h$, and the dot notation as before, observing that $\partial_t h_t=\X$ for all $t$.

\subsection{Perturbation of a simple eigenvalue}\label{sect:perturbation}
We will end this section by finding a formula for the derivative of a simple NEV with respect to a boundary perturbation, in the case that $\lambda$ is a simple NEV for $-\Delta+V$ in $\Omega$. The formula is contained in the following proposition.

  \begin{prop}\label{prop:perturbation}
    Let $\Omega$ be an open domain of class $\C^3$, and take a smooth open domain $\Omega$ such that $\Vcal\subset\Omega$. Let $\lambda_k>\lambda^V$, be a simple eigenvalue for $-\Delta+V$ in $\Omega$, and let $u_k\in H^1(\Omega)$ be its associated normalized eigenfunction, where $\lambda^V$ is the lower bound required in Theorem \ref{thm:direct}. 
    
    Fix a smooth vector field $\X:\Rd\to\Rd$ supported in $\R^n\setminus\Vcal$, and let $h_t(x)=x+t\X(x)$. There exists a differentiable function $t\mapsto\lambda(t)$ for $t$ small enough, such that $\lambda(t)$ is an eigenvalue for $-\Delta+V$ in $\Omega_t=h_t(\Omega)$ and $\lambda(0)=\lambda_k$. 
    
    Furthermore, its derivative at $t=0$ is given by
    \begin{equation}\label{eq:nevpert}
      \dot \lambda(0)=\int_{\partial \Omega} \left(|\nabla_{\partial\Omega}u_k|^2-\lambda_k u_k^2\right)\X\cdot \nu.
    \end{equation}
  \end{prop}
  \begin{proof}
    Consider the map
  \begin{align*}
    F:\left(H^1(\Omega)\cap H^2(\Oprime)\right)&\times I_V\times\mathbb{D}_{\X}\to \left(H^{-1}_0(\Omega)\cap L^2(\Oprime)\right)\times H^{1/2}(\partial \Omega) \times \R\\
    (u,\lambda,h)\mapsto &\left(h^*(\Delta+\lambda-V)h_*u, \,h^*\restr{\nu_{h(\Omega)}\cdot\nabla (h_*u)}{\partial\Omega},\,\int_{\Omega}u^2\right), 
  \end{align*}
  where \[\mathbb{D}_\X= \{h_t\in\textrm{Diff}^3(\Omega):h_t=i_\Omega+t\X,\, t\in(-T,T)\},\] $T^{-1}=\norm{\X}_{\C^1}$ and $I_V=(\lambda^V,+\infty)$. Note that $\mathbb{D}_\X$ is a $1$-dimensional manifold, isomorphic to $(-T,T)$, and that $F$ is $\mathcal{C}^2$.

  Now, we have that $F(u_k,\lambda_k,\textrm{id}_\Omega)=(0,0,1)$ and, whenever $F(u,\lambda,h)=(0,0,1),$ $\lambda$ is a Neumann eigenvalue for $-\Delta+V$ in $h(\Omega)$, with associated normalized eigenfunction $v\defeq h_*u$. Also, since $\lambda_k$ is a simple eigenvalue, we can show that
  \begin{align*}
    \frac{\partial F}{\partial (u,\lambda)}(u_k,\lambda_k,\textrm{id}_\Omega):&\\\left(H^1(\Omega)\cap H^2(\Oprime)\right)\times\R&\to \left(H^{-1}_0(\Omega)\cap L^2(\Oprime)\right)\times H^{1/2}(\partial \Omega) \times \R\\
   (\dot u,\dot \lambda)&\mapsto \left((\Delta+\lambda_k-V)\dot u+\dot\lambda u_k, \,\partial_\nu\dot u,\,2\int_{\Omega}u_k\dot u\right), 
  \end{align*} is an isomorphism. Here we have used \eqref{eq:derivativeL} and \eqref{eq:derivativeN} to compute the derivative.
  
  To see this, observe that the problem
  \begin{equation}\label{eq:udot}
    \begin{cases}
        \left(\Delta +\lambda_k-V\right)\dot u + \dot\lambda u_k = f_1 & \textrm{in }\Omega,\\
        \partial_\nu \dot u= f_2 & \textrm{in } \partial\Omega,\\
        \int_\Omega u_k\, \dot u = \alpha. &
      \end{cases}
  \end{equation}
  will have a solution if and only if
  \[\int_{\partial\Omega}f_2\, u_k= \int_{\Omega} (f_1-\dot \lambda u_k) \,u_k = \int_{\Omega} f_1 u_k - \dot \lambda,\]
  where we have used corollary \ref{cor:np2}, the fact that $\lambda_k$ is a simple eigenvalue with associated eigenfunction $u_k$, and $\int_\Omega u_k^2=1$. This gives a condition that uniquely determines $\dot \lambda$. Now, if we choose a solution $v$ to \eqref{eq:udot}, it will hold that any other solution will be of the form $\dot u = v + su_k$, with $s\in \R$. Applying the last equation in the system \eqref{eq:udot} gives
  \[\int_\Omega u_k\,v+s\int_\Omega u_k^2 = \alpha,\]
  so that $s$ is uniquely determined as \[s=\alpha-\int_\Omega u_k\,v,\]
  and therefore the system above has a unique solution, which proves that\\ $\frac{\partial F}{\partial (u,\lambda)}(u_k,\lambda_k,\textrm{id}_\Omega)$ is an isomorphism, since it is continuous.

  Then, the implicit function theorem says the equation $F(u,\lambda,h)=(0,0,1)$ has a solution $(u_h,\lambda_h)$ for every $h$ in a neighborhood of $i_\Omega$ and that $h\mapsto (u_h,\lambda_h)$ is $\mathcal{C}^1$. Remember that $h_t(x)=x+t\X(x)$, and parametrize $u(t,\centerdot)=u_{h_t}$ and $\lambda(t)=\lambda_{h_t}$. 
  For $t$ close to $0$, we have that $u(t,\centerdot)$ is a normalized eigenfunction for $-\Delta+V$ in $\Omega_t=h_t(\Omega)$ with eigenvalue $\lambda(t)$. 
  
  Now, remember the definition of the anticonvective derivative \eqref{eq:anti-convective}. Although we omit it in the writing, note that $u$, $\lambda$ and $h$ are functions on $t$. We have that
  \[D_t \left[h^*(\Delta+\lambda-V)h_*u\right]=\dot \lambda  u+h^*(\Delta+\lambda-V)h_*\left(D_tu\right),\]
  so at t=0, with $\dot u = \restr{\partial_t u}{t=0}$, we have
  \begin{equation}\label{eq:lambdadot1}
   0=\dot \lambda(0) u_k + (\Delta+ \lambda_k-V)(\dot u \X\cdot \nabla u_k)\quad \textrm{in } \Omega.
   \end{equation}
  Also, we have that
  \[
    \begin{aligned}
      D_t\left[ h^*\left(\nu_{\Omega_t}\cdot\nabla (h_*u)\right)\right]=& h^*\left(\nu_{\Omega_t}\cdot\nabla(h_*D_tu)\right)\\
      &+D_t(h^*\nu_{\Omega_t})\cdot h^*\nabla(h_*u),
    \end{aligned}
    \]
    where $\nu_{\Omega_t}$ is the normal vector to $\partial\Omega_t$ pointing outwards. 
    
    The equality above holds pointwise near $\partial\Omega$ when all the functions involved are smooth, and by continuity it also holds as an equality betweeen continuous curves with values in $H^1(U)$, for $U$ a neighborhood of $\partial\Omega$ in $\R^n$. 
  
  Taking $t\to 0$ above and restricting to $\partial\Omega$, we obtain
  \[-\X\cdot\nabla(\partial_{\nu} u_k)=\partial_\nu(\dot u-\X\cdot\nabla u_k)-(\nabla_{\partial\Omega}\sigma+\sigma\partial_\nu \nu),\]
  where $\sigma=\X\cdot \nu_\Omega$. Since $\partial_\nu u_k=0$ on $\partial\Omega$, we have that $\nabla(\partial_{\nu} u_k)=\partial_\nu(\partial_{\nu}u_k)\nabla u_k$, and therefore
  \begin{align*}
    \partial_\nu(\dot u-\X\cdot\nabla u_k)&=-\sigma\partial_{\nu}(\partial_{\nu} u_k)+(\nabla_{\partial\Omega}\sigma+\sigma\partial_\nu \nu)\nabla u_k\\
    &=\textrm{div}_{\partial\Omega}(\sigma\nabla_{\partial\Omega} u_k)+\sigma\lambda_k u_k,
  \end{align*}
  where we have used that $-\Delta u_k=\lambda_ku_k$ in $\Omega$. If we multiply equation \eqref{eq:lambdadot1} by $u_k$ and integrate over $\Omega$, we obtain 
  \begin{align*}
    0 = & \dot \lambda (0)\int_\Omega u_k^2 + \int_\Omega (\Delta+ \lambda_k-V)(\dot u \X\cdot \nabla u_k)u_k\\
     = & \dot \lambda(0) + \int_{\partial \Omega} u_k\partial_\nu (\dot u- \X\cdot\nabla u_k)\\
      = & \dot \lambda(0) + \int_{\partial \Omega} u_k\textrm{div}_{\partial\Omega}(\sigma\nabla_{\partial\Omega} u_k)+\sigma\lambda_k u_k^2\\
      = & \dot \lambda(0) + \int_{\partial \Omega} \sigma \left(\lambda_k u_k^2-|\nabla_{\partial\Omega}u_k|^2\right),
  \end{align*}
  where we have used that $\int_{\partial \Omega} u_k\textrm{div}_{\partial\Omega}(\sigma\nabla_{\partial\Omega} u_k)=-\int_{\partial\Omega}\sigma|\nabla_{\partial\Omega}u_k|^2$. Finally, we get the desired expression for the derivative of the eigenvalue. \end{proof}

  \section{Simple eigenvalues are generic}\label{sect:nev}

In this section, we will prove that most perturbations produce domains in which Neumann eigenvalues are simple, i.e., that the property of having simple eigenvalues is generic in the appropiate set of perturbations. This will allow us to be in the position to use Proposition \ref{prop:perturbation} to avoid Neumann eigenvalues. To prove it, we will use Henry's genericity theorem which we introduce below. The reader might want to recall some notions on Banach manifolds in Appendix \ref{sect:banach}.

\subsection{Henry's genericity theorem}\label{sect:henry}

Henry gave a generalization of Sard-Smale's theorem in \cite{Henry2005}. We use it in Section \ref{sect:genericity} to prove that ``most perturbations'' on the domain produce simple NEVs. 
%In Section \ref{sect:genericity}, we will use a generalization of Sard-Smale theorem proved by Henry in \cite{Henry2005}. 

Sard proved in 1942 \cite{Sard1942} that the set of critical values of a $\C^k$ map between differentiable manifolds of finite dimension has measure zero, when $k$ is big enough. Smale generalized this result in 1965 \cite{Smale1965} to maps between more general Banach manifolds. In particular, he proved that the set of critical values of a $\C^k$ Fredholm map between Banach manifolds is meager. 

The original version of Smale-Sard theorem can be found in \cite{Smale1965}. Other classical references for the genericity of eigenvalues are \cite{Micheletti1973,uhlenbeck1972,Uhlenbeck1976}. We also used the references \cite{pereiradan2022,Pereira} for better understanding of Henry's work on boundary perturbation.
%He used the notion of transversality, which can be seen as the \textit{opposite} of tangency. \textcolor{red}{Either elaborate or remove this last sentence.}

Henry's generalization can be stated as follows:

  \begin{theorem}[\cite{Henry2005}, Thm. 5.4]\label{thm:trans}
    Let $X,Y,Z$ be Banach manifolds of class $\C^k$, for some $k\in\N$. Let also $A\subset X\times Y$ be open, and $f:A\to Z$ be a map of class $\C^k$. Take a point $\zeta\in Z$ and suppose that the set $f^{-1}(\zeta)$ is Lindel\"of and that, for all $(x,y)\in f^{-1}(\zeta)$,
    \begin{enumerate}
      \item $\partial_x f(x,y):T_xX\to T_\zeta Z$ is semi-Fredholm with index $<k$, and
      \item $df(x,y)=(\partial_x f, \partial_y f):T_xX\times T_yY\to T_\zeta Z$ is surjective.
    \end{enumerate}
    Then, the set $Y_{\textup{crit}}=\{y:\zeta \textrm{ is a critical value of } f(\centerdot\,,y):A_y\to Z\}$ is a meager set in $Y$. Here, $A_y=\{x\in X:(x,y)\in A\}$.
  \end{theorem}

Let's clarify some notions that appear in the theorem. 
\begin{definition}
  A subset $E$ of a topological space $S$ is said to be \textbf{meager} in the sense of Baire (or of first category) if it can be written as a countable union of sets that are nowhere dense in $S$, this is, \[E=\bigcup\limits_{i=1}^\infty E_i, \quad\textrm{ with } \textrm{int}\big(\overline{E_i}\big)=\emptyset.\]  
\end{definition}

The notion of Lindelöf set is a generalization of the notion of compactness:
\begin{definition}
  We say that a topological space has the Lindelöf property if every open cover has a countable subcover.
\end{definition}  

In particular, any separable metric space is Lindelöf, as well as any of its closed subsets.
We should also recall the notion of semi-Fredholm operator:
\begin{definition}
  An operator between Banach spaces $T:X\to Y$ is called \textbf{semi-Fredholm} of index $k\in[-\infty,\infty]$ if $TX$ is closed, and either its kernel is finite-dimensional and/or range has finite codimension, and
\[k\defeq\dim\ker\left(T\right)-\textup{codim}\,\textup{Ran}\left(T\right).\]
\end{definition}

\subsection{Generic simplicity of eigenvalues}\label{sect:genericity}

  Now, we have to carefully define the appropiate spaces of perturbations. In \cite{Henry2005}, Henry chooses the space $\textup{Diff}^3(\Omega)$, the set of embeddings $h:\Omega\to\Rd$ of class $\C^3$ with $\C^3$ inverse.  
  
  However, we have to note that, to avoid regularity problems, the potentials $V$ should not be affected by the deformation. We would like to restrict ourselves to perturbations that equal the identity around $\supp V$. To achieve this we fix from now a cut-off $\chi_1$ in $\C^\infty(\R^n;[0,1])$ such that $\chi_1=0$ in $\supp V$ and $\chi_1=1$ in a neighborhood of $\partial\Omega$. 

  Also, we would like that an open set of the observation set $\Sigma$ is contained in the perturbed domain $h(\Omega)$. Therefore, we will later fix another cut-off $\chi_2$ in $\C^\infty(\R^n;[0,1])$ such that $\chi_2=0$ on a neighborhood of some point $x_0\in\Sigma$.
  
  Then, we may call $\eta=\chi_1\chi_2$. We will consider vector fields $\X\in\C^3\big(\overline\Omega;\R^n\big)$ such that $\norm{\eta\X}_{\C^1}<1,$
  so that any perturbation of the form $h(x)=x+\eta(x)\X(x)$ belongs to $\textup{Diff}^3(\Omega)$. We denote the set of such perturbations as $\mathbb{D}$ and, if we set $M=\norm{\eta}_{\C^3}$, it can be written as follows
  \begin{equation}\label{eq:mathbbd}
      \mathbb{D}=\left\{h:\Omega\to\R^n: h=i_\Omega+\eta\X,\,\X\in\C^3\big(\overline\Omega;\R^n\big),\,\norm{\X}_{\C^3}<M^{-1}\right\},
\end{equation}
where $i_\Omega$ is the identity in $\Omega$.

\begin{figure}[H]
  \centering
  \includegraphics[width=0.7\textwidth]{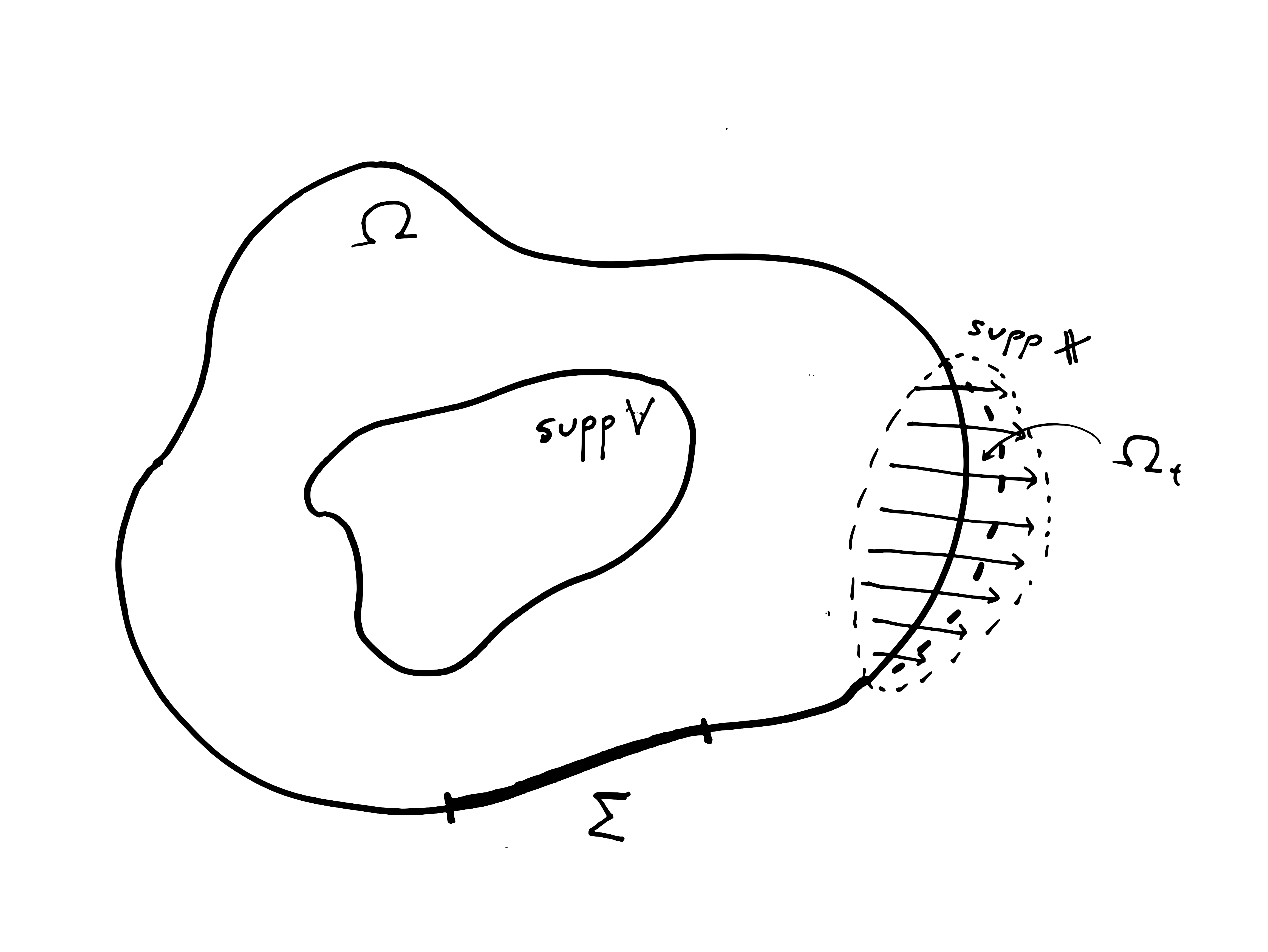}
  \caption{Schematic representation of an admissible deformation}
\end{figure}
    
Note that $\mathbb{D}$ is a Banach manifold of class $\C^3$, since it is isomorphic to an open ball in the Banach space $\C^3\big(\overline\Omega;\R^n\big)$. 

    It will be useful to characterize the tangent space of $\mathbb{D}$ at the identity $i_\Omega(x)=x$. Observe that any curve in $t\mapsto h_t \in \mathbb{D}$ such that $h_0(x)=x$ can be written as $h_t(x)=x+t\eta(x)\X(x)$ for some vector field $\X\in\C^3\big(\overline\Omega;\R^n\big)$, and $t\in(-T,T)$, where \[T=\norm{\eta(x)\X(x)}_{\C^1}^{-1}.\] Differentiating with respect to $t$ at $t=0$ gives the shape of the tangent space at $i_\Omega$: 
    \[T_{i_\Omega}\mathbb{D}=\left\{\dot h = \eta\X:\X\in\C^3\big(\overline\Omega;\R^n\big)\right\}.\]

We are now ready to state the main result of this section:
%    On the other hand, our ultimate goal in the solution of the inverse problem will be to still need part of the observation sets $\Sigma_1$ and $\Sigma_2$ to be contained in the perturbed domain $h(\Omega)$. Therefore, we will need to consider perturbations that are the identity in a neighborhood of a point in $\Sigma_1$ and another in $\Sigma_2$. Indeed, set $\Sigma_1,\Sigma_2\subset\partial\Omega$ relatively open, and define
%    \[\Eta=\{\eta\in \C^\infty(\R^n;[0,1]):\exists x_j\in \Sigma_j \textup{ such that } \eta=0 \textup{ in a nhood. of } x_j,\,j=1,2\}.\]
%    
%    Therefore, the space we are interested in is
%    \[\mathfrak{D}=\{h_t=h'(t,\centerdot)\eta:h'(t,\centerdot)\in\mathbb{D},\,\eta\in\Eta\}.\]
  \begin{prop}\label{prop:genericity}
     Let $\Omega$ be a bounded open domain of class $\C^3$ such that $\supp V\subset \Omega$, anf fix $\chi_1$ as discussed above. Let $\lambda_k>\lambda^V$ be a Neumann eigenvalue for $-\Delta+V$ in $\Omega$.

    Then, there exists a cut-off $\chi_2\in\C^\infty(\R^n;[0,1])$ such that $\chi_2=0$ on a neighborhood of some point $x_0\in\Sigma$, such that the set
    \[\{h\in \mathbb{D} : \textup{$\lambda_k$ is a multiple eigenvalue for} -\Delta+V \textup{ in } h(\Omega)\}\]
    is meager in $\mathbb{D}$, as defined in \eqref{eq:mathbbd} with $\eta=\chi_1\chi_2$.
  \end{prop}
  We will divide the proof of this proposition in a series of lemmas, to hopefully make it more clear. We again denote by $\Vcal$ the support of the potential $V$.
  \begin{lemma}\label{lem:gen1}
    Let $\Omega$ and $\Vcal$ be as above, and let $\lambda_k>\lambda^V$ be a multiple eigenvalue for $-\Delta+V$ in $h(\Omega)$ for some $h\in\mathbb{D}$. Then, $(0,0)$ is a singular point for the map 
    \begin{align*}
      F_h:\left(H^1(\Omega)\cap H^2(\Oprime)\setminus\{0\}\right)&\times I_V\to \left(H^{-1}_0(\Omega)\cap L^2(\Oprime)\right)\times H^{1/2}(\partial \Omega) \\
      (u,\lambda)\mapsto& \left(h^*(\Delta+\lambda-V)h_*u, \,h^*\restr{\nu_{h(\Omega)}\cdot\nabla (h_*u)}{\partial\Omega}\right), 
    \end{align*}
    where $I_V=(\lambda^V,+\infty)$.
  \end{lemma}

  \begin{proof}
    By performing a ``change of origin'', we can reduce to the case of $h=i_\Omega$. Then, suppose that $\lambda_k$ is a multiple Neumann eigenvalue for $-\Delta+V$ in $\Omega$, we can show that the map
    \begin{align*}
      F_{i_\Omega}:\left(H^1(\Omega)\cap H^2(\Oprime)\setminus\{0\}\right)\times I_V&\to \left(H^{-1}_0(\Omega)\cap   L^2(\Oprime)\right)\times H^{1/2}(\partial \Omega) \\
      (u,\lambda)&\mapsto \left((\Delta+\lambda-V)u, \,\restr{\partial_\nu u}{\partial\Omega}\right), 
    \end{align*}
   % at $(u,\lambda)\in\left(H^1(\Omega)\cap H^2(\Oprime)\setminus\{0\}\right)\times\R$ is
    has $(0,0)$ as a singular value. For this take $u_k$ to be any eigenfunction associated to $\lambda_k$, and note that $(u_k,\lambda_k)\in F_{i_\Omega}^{-1}(0,0)$.
    
    The differential of $F_{i_\Omega}$ at $(u_k,\lambda_k)$
    \begin{align*}
      dF_{i_\Omega}(u_k,\lambda_k):\left(H^1(\Omega)\cap H^2(\Oprime)\right)\times\R&\to \left(H^{-1}_0(\Omega)\cap L^2(\Oprime)\right)\times H^{1/2}(\partial \Omega) \\
      (\dot u,\dot \lambda)&\mapsto \left((\Delta+\lambda_k-V)\dot u +\dot\lambda u_k, \,\restr{\partial_\nu \dot u}{\partial\Omega}\right),
    \end{align*}
    will be surjective only if the problem 
    \begin{equation}\label{eq:np3}
      \begin{cases}
          \left(\Delta +\lambda_k-V\right)\dot u + \dot\lambda u_k = f_1 & \textrm{in }\Omega,\\
          \partial_\nu \dot u= f_2 & \textrm{in } \partial\Omega
        \end{cases}
      \end{equation}
      has a solution for every $(f_1,f_2)\in \left(H^{-1}_0(\Omega)\cap L^2(\Oprime)\right)\times H^{1/2}(\partial \Omega)$. By corollary \ref{cor:np2}, this would mean that
      \[ \int_{\partial \Omega} f_2\, u = \int_{\Omega} (f_1-\dot \lambda u_k)\, u\]
  for all $u\in H^1(\Omega)$ solution to the homogeneous Neumann problem \eqref{eq:neuhom}. Since we have supposed $\lambda_k$ to be a multiple eigenvalue, there must exist a non-zero eigenfunction $u_k'$ orthogonal to $u_k$. This means that
  \[ \int_{\partial \Omega} f_2\, u_k' = \int_{\Omega} f_1\, u_k'\]
  for all $(f_1,f_2)\in \left(H^{-1}_0(\Omega)\cap L^2(\Oprime)\right)\times H^{1/2}(\partial \Omega)$ and, in particular, taking $f_2=0$ we get
  \[ \int_{\Omega} f_1\, u_k'=0\]
  for all $f_1\in L^2(\Omega)$, which means that $u'_k=0$. This is a contradiction with $u_k'$ being an eigenfunction.
  \end{proof}

  \begin{lemma}\label{lem:gen2}
    Let $F_h$ be as in Lemma \ref{lem:gen1}. If $(u_k,\lambda_k)\in F_h^{-1}(0,0)$, then $dF_h(u_k,\lambda_k)$ is Fredholm with index 1, this is,
    \[\dim\ker\left(dF_h(u_k,\lambda_k)\right)-\textup{codim}\,\textup{Ran}\left(dF_h(u_k,\lambda_k)\right)=1.\]
  \end{lemma}
  \begin{proof}
    Perform again a ``change of origin'' and observe that the kernel of \linebreak $dF_{i_\Omega}(u_k,\lambda_k)$ is the space of solutions $(\dot u,\dot \lambda)$ to the problem
    \begin{equation*}
      \begin{cases}
          \left(\Delta +\lambda_k-V\right)\dot u + \dot\lambda u_k = 0 & \textrm{in }\Omega,\\
          \partial_\nu \dot u= 0 & \textrm{in } \partial\Omega,
        \end{cases}
      \end{equation*}
      to which corollary \ref{cor:np2} imposes \[\int_\Omega \dot \lambda u_k^2 = 0,\] and therefore $\dot \lambda=0$. This means that the kernel of $dF_{i_\Omega}(u_k,\lambda_k)$ is equal to the eigenspace of $\lambda_k$, whose dimension can be denoted by $m_{\lambda_k}.$ On the other hand, the image of $dF_{i_\Omega}(u_k,\lambda_k)$ is given as those $(f_1,f_2)$ such that problem \eqref{eq:np3} has a solution, and we have previously seen that this is equivalent to the condition
      \[ \int_{\partial \Omega} f_2\, u = \int_{\Omega} f_1\, u\]
      for every $u$ Neumann eigenvalue orthogonal to $u_k$. Therefore, its codimension is given by $m_{\lambda_k}-1$, and thus the index of $dF_{i_\Omega}(u_k,\lambda_k)$ is 1.
  \end{proof}

  \begin{lemma}\label{lem:gen3}
    Let $\Omega$ be as in Proposition \ref{prop:genericity}, let $\lambda_k>\lambda^V$ be a Neumann eigenvalue for $-\Delta-V$ in $\Omega$, and let $u_k$ be an associated eigenfunction. 
    
   If $\phi\in \C^2(\partial\Omega)$ is such that 
   \begin{equation*}
    \nabla_{\partial\Omega}\phi\cdot \nabla_{\partial\Omega} u_k + \phi\,\Delta_{\partial\Omega}u_k+\phi\,\partial^2_\nu u_k=0 \qquad\textup{in }\partial\Omega.
  \end{equation*}
  Then, $\phi=0$.
  \end{lemma}
\begin{proof}
  First, observe that, by theorem 11.1.1 in \cite{hormanderlinear1963} and theorem 6.3.2.1 in \cite{Grisvard2011}, we have that $u_k\in \C^2(\partial\Omega)$. Now, we have the following formula for the Laplace-Beltrami operator:
    \[\Delta_{\partial\Omega}u_k=\Delta u_k-H \partial_\nu u_k-\partial_\nu^2 u_k,\] where $H=\nabla\cdot\nu$ is the mean curvature of $\partial\Omega$. Since near $\partial\Omega$ we have that $\Delta u_k+\lambda_k u_k=0$ and $\partial_\nu u=0$, we arrive to the expression 
    \begin{equation}\label{eq:contrapsi3}
      \nabla_{\partial\Omega}\phi\cdot \nabla_{\partial\Omega} u_k=\lambda_k\phi\,u_k\qquad\textup{in }\partial\Omega.
    \end{equation}

    Now, since $\partial\Omega$ is compact and $u_k,\phi\in\C^2\big(\partial\Omega)$, we can find $x_0\in \textup{argmax}_{\partial\Omega} \lvert u_k\phi\rvert$. Then, we can define a curve $t\mapsto x(t)$ in $\partial\Omega$ by $x(0)=x_0$ and $x'(t)=\nabla_{\partial\Omega}u_k(x(t))$ in some interval $t\in(-M,M)$. We have by the chain rule that \[\frac{\dd}{\dd t} u_k(x(t))=\nabla_{\partial\Omega} u_k(x(t))\cdot \nabla_{\partial\Omega}u_k(x(t))=|\nabla_{\partial\Omega}u_k(x(t))|^2\geq 0,\] and that \[\frac{\dd}{\dd t} \phi(x(t)) = \nabla\phi(x(t)) \cdot \nabla_{\partial\Omega}u_k(x(t))= \lambda_k\phi(x(t))u_k(x(t)),\] where we have also used \eqref{eq:contrapsi3}. This means that \[\phi(x(t)) = \phi(x_0)\exp\left(\lambda_k\int_0^t u_k(x(s))\dd s\right).\] Now, if $u_k(x_0)>0$, we have that $\lvert\phi(x(t))\rvert>\lvert\phi(x_0)\rvert$ for $t>0$, and therefore $\lvert u_k(x(t))\phi(x(t))\rvert > \lvert u_k(x_0)\phi(x_0)\rvert$ for $t>0$. This is a contradiction with the fact that $\lvert u_k\phi\rvert$ has a maximum at $x_0$.  

    On the other hand, if $u_k(x_0)<0$, we have that $\lvert\phi(x(t))\rvert>\lvert\phi(x_0)\rvert$ for \linebreak$t<0$, and therefore $\lvert u_k(x(t))\phi(x(t))\rvert > \lvert u_k(x_0)\phi(x_0)\rvert$ for $t<0$. This is again a contradiction. 

    Therefore, we must have that $u_k\phi=0$ in $\partial\Omega$. However, if $\phi$ is not identically zero, we can find an open subset of $\partial\Omega$ in which $u_k=0$. By the unique continuation property, Corollary \ref{prop:ucp}, we will have that $u_k=0$ in $\Omega$, which is a contradiction with the fact that $u_k$ is a Neumann eigenvalue. Therefore, $\phi=0$.
\end{proof}

  We are now ready to prove the proposition.
  \begin{proof}[\textbf{Proof of Proposition \ref{prop:genericity}}:]
    Define the following function
    \begin{align*}
      F:\left(H^1(\Omega)\cap H^2(\Oprime)\setminus\{0\}\right)&\times I_V\times\mathbb{D}\to \left(H^{-1}_0(\Omega)\cap L^2(\Oprime)\right)\times H^{1/2}(\partial \Omega) \\
      (u,\lambda,h)\mapsto& \left(h^*(\Delta+\lambda-V)h_*u, \,h^*\restr{\nu_{h(\Omega)}\cdot\nabla (h_*u)}{\partial\Omega}\right), 
    \end{align*}
    where again $I_V=(\lambda^V,+\infty)$.
    
  By Lemma \ref{lem:gen1}, we just have to show that the set of such $h\in\mathbb{D}$ for which $(0,0)$ is a critical value of $F_h=F(\centerdot\,,\centerdot\,,h)$ is meager in $\mathbb{D}$, and for this we use Theorem \ref{thm:trans}.  
  
  We'll check that the hypotheses of the theorem are fullfilled. Indeed, by Lemma \ref{lem:gen2}, we know that, if $(u_k,\lambda_k)\in F_h^{-1}(0,0)$, then $dF_h(u_k,\lambda_k)$ is Fredholm with index 1. Note also that $F^{-1}(0,0)$ is Lindel\"of, since it is a subset of a separable space.

    Therefore, the only thing left to show is that, whenever $(u_k,\lambda_k,h_0)\in F^{-1}(0,0)$, this is, whenever $\lambda_k$ is a Neumann eigenvalue for $-\Delta+V$ in $h_0(\Omega)$ with $u_k$ an associated eigenfunction, the differential of $F$ at $(u_k,\lambda_k,h_0)$ is surjective. Again, we can ``change the origin'' to suppose $h_0=i_\Omega$. The differential of $F$ at $(u_k,\lambda_k,i_\Omega)$ can be computed using \eqref{eq:derivativeL} and \eqref{eq:derivativeN}, and it's given by
    \begin{align*}
      dF(u_k,\lambda_k,i_\Omega):&\\\big(H^1(\Omega)&\cap H^2(\Oprime)\big)\times\R\times T_{i_\Omega}\mathbb{D}\to \left(H^{-1}_0(\Omega)\cap L^2(\Oprime)\right)\times H^{1/2}(\partial \Omega) \\
      &\hspace{2.85cm}(\dot u,\dot \lambda,\dot h)\mapsto\big( G_1, \,G_2\big),
    \end{align*}
    where \[G_1=(\Delta+\lambda_k-V)(\dot u-\dot h\cdot \nabla u_k) +\dot\lambda u_k,\] and \[G_2=\partial_\nu(\dot u-\dot h\cdot \nabla u_k)-\nabla_{\partial\Omega}(\dot h\cdot \nu)\cdot\nabla u_k + \dot h\cdot\nabla (\partial_\nu u_k).\] 

    We will proceed by contradiction. Suppose that $dF(u_k,\lambda_k,i_\Omega)$ is not surjective, then there must exist $(\psi,\theta)\in \left(H^{-1}_0(\Omega)\cap L^2(\Oprime)\right)\times H^{1/2}(\partial \Omega)$ orthogonal to the image of $dF(u_k,\lambda_k,i_\Omega)$, and $(\psi,\theta)\neq(0,0)$. This means that for all $(\dot u,\dot \lambda,\dot h)\in \big(H^1(\Omega)\cap H^2(\Oprime)\big)\times\R\times T_{i_\Omega}\mathbb{D}$ we must have 
    \[\int_\Omega\psi \, G_1 + \int_{\partial\Omega} \theta\,G_2 = 0.\]
    If we take $\dot h = \dot \lambda = 0$, we have that 
    \begin{equation}\label{eq:psitheta}
       \int_\Omega \psi \,(\Delta+\lambda_k-V)\dot u+\int_{\partial\Omega}\theta\,\partial_\nu \dot u= 0,
    \end{equation}
    for all $\dot u\in H^1(\Omega)\cap H^2(\Oprime)$. In particular,    \[ \int_\Omega \psi \,(\Delta+\lambda_k-V)\dot u= 0,\] 
    for all $\dot u \in C_0^\infty(\Omega)$. This means $(\Delta+\lambda_k-V)\psi=0$ in $\Omega$ in the distributional sense, which implies that $\psi\in H^1(\Omega)\cap H^2(\Oprime)$. Integrating by parts in \eqref{eq:psitheta} we obtain
    \begin{equation}\label{eq:psithetabis} 
      \int_{\partial\Omega}\psi\,\partial_\nu\dot u-\int_{\partial\Omega}\partial_\nu\psi\,\dot u+\int_{\partial\Omega}\theta\,\partial_\nu\dot u=0,
    \end{equation}
     and choosing any $\dot u$ such that $\partial_\nu \dot u=0$ in $\partial\Omega$ yields
    \[\int_{\partial\Omega}\partial_\nu \psi \,\dot u=0 \] 
    for all $\dot u \in H^{3/2}(\partial\Omega)$, and by density for all $\dot u\in L^2(\partial\Omega)$. This means that $\partial_\nu\psi=0$ and therefore $\psi$ solves
    \begin{equation}\label{eq:psinev} \begin{cases}
        \left(\Delta+\lambda_k-V\right)\psi = 0 & \textrm{in }\Omega,\\
        \partial_\nu \psi = 0 & \textrm{in }\partial\Omega,
      \end{cases}\end{equation}
    Therefore, from \eqref{eq:psithetabis}, we have that
    \begin{equation}\label{eq:psitheta2} 
      \int_{\partial\Omega}\psi\,\partial_\nu\dot u=-\int_{\partial\Omega}\theta\,\partial_\nu\dot u,
    \end{equation}
    for all $\partial_\nu\dot u\in H^{1/2}(\partial\Omega)$ and thus $\restr{\psi}{\partial\Omega}=-\theta$.
    
    If we take $\dot h = \dot u = 0$, then we have that
    \begin{equation*}
      \int_\Omega u_k\,\psi = 0.
    \end{equation*}

    Finally, if $\dot u = \dot \lambda = 0$, we obtain by \eqref{eq:psitheta2},
    \begin{align*}
       0=&\int_\Omega \psi\,(\Delta+\lambda_k-V)(\dot h\cdot\nabla u_k) \\ &- \int_{\partial\Omega}\psi\left[\partial_\nu(\dot h\cdot \nabla u_k)+\nabla_{\partial\Omega}(\dot h\cdot \nu)\cdot\nabla u_k - \dot h\cdot\nabla (\partial_\nu u_k)\right] = 0.
    \end{align*}

    Note that $u,\psi\in \C^\infty(\Oprime)$ by theorem 11.1.1 in \cite{hormanderlinear1963}, and that the domains and perturbations are $\C^3$. Therefore, by theorem 6.3.2.1 in \cite{Grisvard2011}, we have that $u_k$ and $\psi$ are in $\C^{2,\alpha}\big(\overline{\Omega\setminus\Vcal}\big)$, $0<\alpha<1$. Integrating by parts we see that
    \[\int_\Omega \psi\,(\Delta+\lambda_k-V)(\dot h\cdot\nabla u_k) - \int_{\partial\Omega}\psi\,\partial_\nu(\dot h\cdot \nabla u_k)=0,\]
    and therefore calling $\sigma = \dot h \cdot \nu$, we obtain
    \[\int_{\partial\Omega}\psi\left[\nabla_{\partial\Omega}\sigma\cdot\nabla u_k - \sigma\,\partial^2_\nu u_k\right]=0,\]
    where we have used that $\nabla(\partial_\nu u_k)=\nu\,\partial^2_\nu u_k$ since $\restr{\partial_\nu u_k}{\partial\Omega}=0$ and thus \linebreak $\nabla_{\partial\Omega}(\partial_\nu u_k)=0$. 
    
    Observe now that, since $\partial_\nu u_k=0$, then $\restr{\nabla u_k}{\partial\Omega}=\nabla_{\partial\Omega} u_k$, and therefore 
    \[\psi\,\nabla_\Omega \sigma\cdot \nabla u_k = \textup{div}_{\partial\Omega}(\sigma\psi\nabla_{\partial\Omega}u_k)-\sigma(\nabla_{\partial\Omega}\psi\cdot \nabla_{\partial\Omega} u_k + \psi\Delta_{\partial\Omega}u_k).\]Then, we have that 
    %\begin{equation}\label{eq:contrapsi2}
    %\begin{aligned}
    %  &\int_{\partial\Omega}\sigma\left[\nabla_{\partial\Omega}\psi\cdot \nabla_{\partial\Omega} u_k + \psi\,\Delta_{\partial\Omega}u_k+\psi\,\partial^2_\nu u_k\right]\\
    %  &=\int_{\partial\Omega}\textup{div}_{\partial\Omega}(\sigma\psi\nabla_{\partial\Omega}u_k)=0\quad \textup{for all } \sigma=\dot h\cdot\nu, \;\dot h\in T_{i_\Omega}\mathbb{D}
    %\end{aligned}
    %\end{equation}
    \begin{equation}\label{eq:contrapsi2}
      \begin{aligned}
        &\int_{\partial\Omega}\sigma\left[\nabla_{\partial\Omega}\psi\cdot \nabla_{\partial\Omega} u_k + \psi\,\Delta_{\partial\Omega}u_k+\psi\,\partial^2_\nu u_k\right]\\
        &=\int_{\partial\Omega}\textup{div}_{\partial\Omega}(\sigma\psi\nabla_{\partial\Omega}u_k)=0
      \end{aligned}
      \end{equation}
    for all $\sigma=\dot h\cdot\nu$ with $\dot h\in T_{i_\Omega}\mathbb{D}$. Here we have used the fact that $\nabla_{\partial\Omega}u_k\cdot\nu=0$ along with the divergence theorem, this is, for any vector field $\mathbf X\perp\nu$ on $\partial\Omega$ we have that
    \[\int_{\partial\Omega}\textup{div}_{\partial\Omega}(\mathbf X)=0.\]

    %\[\frac{\partial n}{\partial n}= n\cdot \nabla n = 0,\quad\textup{since }n\perp\nabla n \]
    
    Note that by assumption $\psi\neq 0$. Indeed, if $\psi=0$, then $\theta=\restr{\psi}{\partial\Omega}=0$, and we assumed $(\psi,\theta)\neq(0,0)$. Hence, we have by Lemma \ref{lem:gen3} that the term
    \begin{equation}\label{eq:theterm}
    \nabla_{\partial\Omega}\psi\cdot \nabla_{\partial\Omega} u_k + \psi\,\Delta_{\partial\Omega}u_k+\psi\,\partial^2_\nu u_k
    \end{equation}
    can't be zero in the whole of $\partial\Omega$. Then, choose in the definition of $\mathbb{D}$ a cut-off $\chi_2\in\C^\infty(\R^n;[0,1])$ such that $\chi_2\neq0$ in points in which \eqref{eq:theterm} doesn't vanish for any possible $\psi$, and at the same time $\chi_2=0$ is a neighborhood of a point $x_0\in\Sigma$. This is possible because, as seen in \eqref{eq:psinev}, $\psi$ is a Neumann eigenvalue associated to $\lambda_k$, of which there can only be finitely many linearly independent ones, by corollary \ref{cor:finitenev}. 

    Finally, choose an open set $U\in\partial\Omega$ in which \[\left|\nabla_{\partial\Omega}\psi\cdot \nabla_{\partial\Omega} u_k + \psi\,\Delta_{\partial\Omega}u_k+\psi\,\partial^2_\nu u_k\right|>0,\] and take $\dot h = \X \eta$ for some $\X$ supported in $U$ such that $\X\cdot\nu\geq0$. Then, 
    \[\int_{\partial\Omega}\dot h \cdot \nu\left[\nabla_{\partial\Omega}\psi\cdot \nabla_{\partial\Omega} u_k + \psi\,\Delta_{\partial\Omega}u_k+\psi\,\partial^2_\nu u_k\right]\neq 0,\]
    which contradicts \eqref{eq:contrapsi2}.  

    Therefore, $dF(u_k,\lambda_k,i_\Omega)$ is surjective, and we are in position to use Theorem \ref{thm:trans}. This finishes the proof of the proposition.
     
    %We know, by unique continuation, Proposition \ref{prop:ucp}, that neither $u_k$ nor $\psi$ can vanish on a open subset of $\partial\Omega$. Thus, we can choose $x_0\in\partial\Omega$ such that $u_k(x_0)\psi(x_0)\neq 0$. Define now a curve $t\mapsto x(t)$ in $\partial\Omega$ by $x(0)=x_0$ and $\frac{\dd}{\dd t} x(t)=\nabla_{\partial\Omega}u_k(x(t))$. Since $\nabla_{\partial\Omega}u_k$ is $\C^1$ and $\partial\Omega$ is compact, this curve is defined for $t\in\R$, and we have by the chain rule that \[\frac{\dd}{\dd t} u_k(x(t))=\nabla u_k(x(t))\cdot \nabla_{\partial\Omega}u_k(x(t))=|\nabla_{\partial\Omega}u_k(x(t))|^2\geq 0,\] and that \begin{align*}\frac{\dd}{\dd t} \psi(x(t)) = \nabla\psi(x(t)) \cdot \nabla_{\partial\Omega}u_k(x(t))= \lambda\psi(x(t))u_k(x(t)),\end{align*}
    %so that 
    %\[ \psi(x(t)) = \psi(x_0)\exp\left(\lambda\int_0^t u_k(x(s))\dd s\right).\]
    %If we chose $x_0$ such that $u_k(x_0)>0$, we have that $u_k(x(t))\geq u_k(x_0)>0$ for all $t\in\R$, and therefore $|\psi(x(t))|\to\infty$ as $t\to\infty$ as long as $\lambda>0$, which is a contradiction with the fact that $\psi$ is bounded. On the other hand, if $u_k(x_0)<0$, then $|\psi(x(t))|\to\infty$ as $t\to-\infty$, which is again a contradiction. \textcolor{red}{In this part of the argument, we need to be careful with the regularity of the functions. Henry says that the functions are $\C^{2,\alpha}$ as long as the domain is $\C^3$, and then uses a density argument to extend the result for $\C^2$ domains. Check why this is the case, and decide which regularity to impose on domains.}
\end{proof}

\section{Proof of Theorem \ref{thm:nev}}\label{sect:thm}

Let's observe first that eigenvalues in fact form a countable set. To see this, define the unbounded operator $\left(T_N,D(T_N)\right)$ over $L^2(\Omega)$ as $T_Nu=\left(-\Delta+V\right)u$, with domain 
\[ D(T_N)=\{u\in L^2(\Omega): \left(-\Delta+V\right)u\in L^2(\Omega),\; \partial_\nu \restr{u}{\partial \Omega}=0\}.\]
Observe now that a Neumann eigenvalue for $-\Delta+V$ on $\Omega$ will be an eigenvalue for $\left(T_N,D(T_N)\right)$. The domain $D(T_N)$ is a separable Hilbert space, since it is a subspace of $L^2(\Omega)$, and $T_N$ is symmetric over $D(T_N)$, which ensures that the set of its eigenvalues must be countable. 

Indeed, suppose that there is an uncountable set of such eigenvalues. Let $\lambda$ and $\mu$ be any two distinct eigenvalues, and $u$ and $v$ be corresponding distinct eigenfunctions. Then,
\[
  \lambda\braket{u,v}=\braket{T_N u, v}=\braket{u,T_N v}= \mu\braket{u,v},
  \]
and thus $u\perp v$, which contradicts the fact that $D(T_N)$ is separable.
We recall once again the statement of the theorem we want to prove:
\nevthm*
\begin{proof}
To prove the theorem, choose any open domain $\Omega$ such that $\Sigma\subset\partial\Omega$, $\supp V\subset\Omega$. There are various possibilities:
\begin{enumerate}
  \item In case $\lambda$ is \textbf{not} a Neumann eigenvalue for $-\Delta+V$ in $\Omega$, we can use corollary \ref{cor:nev} to solve problem \eqref{eq:neumann}.

  \item In case $\lambda$ is a \textbf{simple} Neumann eigenvalue, call it $\lambda_k$, and let $u_k$ be an associated eigenfunction. We can use Proposition \ref{prop:perturbation} to move this eigenvalue \textit{away} from $\lambda$. Indeed, remember that, if we choose a perturbation of the form $h_t(x)=x+t\X(x)$ for $\X\in\C^3\big(\overline\Omega;\R^n\big)$ supported away from $\supp V$, $t$ small enough, we can obtain a curve of eigenvalues $t\mapsto \mu(t)$ in $\Omega_t=h_t(\Omega)$ such that $\mu(0)=\lambda_k$, and the derivative of this curve at $t=0$ is given by
  \[\dot \mu(0)=\int_{\partial \Omega} \left(|\nabla_{\partial\Omega}u_k|^2-\lambda_k u_k^2\right)\X\cdot\nu.\]

  Observe now that $|\nabla_{\partial\Omega}u_k|^2-\lambda_k u_k^2$ can't vanish in the whole boundary $\partial\Omega$. Indeed, suppose that $|\nabla_{\partial\Omega}u_k|^2=\lambda_k u_k^2$ in $\partial\Omega$, note that $u\in \C^2(\overline\Omega)$, and therefore $u_k$ attains the maximum and the minimum in the compact surface $\partial\Omega$. Let $x_m\in\arg\min_{\partial\Omega}\{u_k\}$ and $x_M\in\arg\max_{\partial\Omega}\{u_k\}$. Then, 
  \[
    \nabla_{\partial\Omega}u_k(x_m)=\nabla_{\partial\Omega}u_k(x_M)=0,
  \]
  which means that $u_k(x_m)=u_k(x_M)=0$, and therefore $u_k\equiv 0$ in $\partial\Omega$, which, by the unique continuation property, Proposition \ref{prop:ucp}, implies that $u_k\equiv 0$ in $\Omega$, which is a contradiction with the fact that $u_k$ is an eigenfunction. Therefore, we can find a point $x_0\in\partial\Omega$ such that $|\nabla_{\partial\Omega}u_k|^2-\lambda_k u_k^2\geq 0$ in a neighborhood of $x_0$. 
  
  If we then choose $\X$ to be supported near $x_0$, and such that $\X\cdot\nu>0$ in its support, we obtain that $\dot{\mu}(0)>0$, and therefore $\lambda\neq\mu(t)$ for some $t\neq0$. Note that the rest of eigenvalues would move but, since they form a discrete set, they won't become equal to $\lambda$ if $t$ is small enough. This means that $\Omega_t$ is a domain in which $\lambda$ is {not} a Neumann eigenvalue of $-\Delta+V$.

  Note that we have to choose $\X$ such that $\X\cdot\nu=0$ in some non-empty open set $\Sigma'\subset\Sigma$, so that the perturbed boundary $\partial \Omega_t$ contains $\Sigma'$.

  \item Finally, suppose $\lambda$ is a \textbf{multiple} eigenvalue for $\Omega$. Then, we only need to choose a perturbation $h\in\mathbb{D}$ as in Section \ref{sect:genericity}, such that $\lambda$ is not a multiple eigenvalue for $-\Delta+V$ in $h(\Omega)$. This can be done thanks to Proposition \ref{prop:genericity}. If $\lambda$ is not a NEV in $h(\Omega)$, we are done, and in case it is, it will be a simple eigenvalue, and we are back to the previous point.
\end{enumerate}
\end{proof}

\appendix

\section{Banach manifolds}\label{sect:banach}
Banach are a generalization of the usual finite dimensional manifolds, where the charts are defined on general Banach spaces instead of on $\R^n$. A general reference for this topic is \cite{Lang1999}. We give here just a brief introduction to suit our needs.

\begin{definition}
    Let $M$ be a set. An \textbf{atlas of class} $\C^k$ ($k\in\N_0$) on $M$ is a collection of pairs $\{(U_i,\varphi_i)\}_{i\in I}$, called \textbf{charts}, where $I$ is an index set, satisfying the following conditions
    \begin{enumerate}
        \item each $U_i$ is a subset of $M$ and $M=\cup_{i\in I} U_i$,
        \item each $\varphi_i$ is a bijection onto an open subset $\varphi(U_i)$ of a Banach space $E_i$, and for any $i,j$ we have that $\varphi_i(U_i\cap U_j)$ is open in $E_i$, and
        \item the map \[\varphi_j\varphi^{-1}_i:\varphi_i(U_i\cap U_j)\to\varphi_j(U_i\cap U_j)\] is an isomorphism of class $\C^k$ for each pair of indices $i,j\in I$.
    \end{enumerate}
\end{definition}

We can easily check that there is a unique topology in $M$ such that each $U_i$ is open and each $\varphi_i$ is an homeomorphism. For $M$ to be Haussdorff, we have to place a separation condition on the covering.

We have not made any assumption on the relationship between the $E_is$, but if they are all equal to the same space $E$, we say that we are in possesion of an \textbf{$E$-atlas}. Furthermore, we can show that, whenever $U_i\cap U_j\neq\emptyset$, $E_i$ and $E_j$ are isomorphic as topological vector spaces, and therefore, we will have an $E$-atlas in any connected component on $M$.

Now suppose that we are given an open subset $U\subset M$ and a topological isomorphism $\varphi$ onto an open subset $\varphi(U)$ of a Banach space $E$. We say that $(U,\varphi)$ is \textbf{compatible} with the atlas $\{(U_i,\varphi_i)\}_{i\in I}$ if \[\varphi\varphi_i^{-1}:\varphi(U_i\cap U)\to\varphi_i(U_i\cap U_j)\] is an isomorphism of class $\C^k$ for all $i\in I$. We say that two atlas are compatible if any chart in one of them is compatible with the other atlas. The relation of compatibility is a equivalence relation on the set of atlases of class $\C^k$.
\begin{definition}
    A \textbf{Banach manifold of class} $\C^k$ is a set $M$ together with an equivalence class of atlases of class $\C^k$.
\end{definition}
If all the spaces $E_i$ in an atlas are isomophormic as topological vector spaces, then we can find an equivalent atlas for which all spaces are all equal, say to $E$. In this case, we say that $M$ is a \textbf{$E$-manifold}, or that $M$ is \textbf{modeled on} $E$. If $E$ is isomorphic to $\R^m$, we say that $M$ is a \textbf{finite dimensional manifold} of dimension $m$. 

We will also need to define tangent vectors to endow the manifold with a differentiable structure. Let $M$ be a $E$-manifold of class $\C^k$, $k\geq 1$, and let $x$ a point in $M$. Consider triples $(U,\varphi,v)$, where $(U,\varphi)$ is an atlas such that $x\in U$, and $v$ is a vector in $E$. We say that two such triples $(U,\varphi,v)$ and $(V,\psi,w)$ are equivalent if 
\[\restr{D(\psi\varphi^{-1})}{\varphi(x)}v = w.\]
\begin{definition}
    An equivalent class of triples $(U,\varphi,v)$ is called a \textbf{tangent vector} at $x$. The set of all tangent vectors at $x$ is the \textbf{tangent space} to $M$ at $x$, and is denoted by $T_xM$.
\end{definition}
Each chart determines a bijection between the tangent space and a Banach manifold, and allows to transport the structure of topological vector space to it. This structure will be independent of the chart.

We can equivalently define tangent vectors as equivalence classes of curves that \textit{pass through $x$ with the same velocity.} Indeed, let $\gamma_1,\gamma_2:(-\varepsilon,\varepsilon)\to M$ be two curves such that $\gamma_1(0)=\gamma_2(0)=x$. We say that $\gamma_1$ and $\gamma_2$ are equivalent at $x$ if
\[(\varphi\gamma_1)'(0)=(\varphi\gamma_2)'(0),\]
for $(U,\varphi)$ a chart containing $x$. Then, we can define the tangent vector $v$ at $x$ as $v=(\varphi\gamma)'(0)$ for any $\gamma$ in the equivalence class.

\bibliographystyle{alpha}
\bibliography{references}
\end{document}